\documentclass[preprint]{elsarticle}
\usepackage{amssymb}
\usepackage{amsfonts}
\usepackage{amsmath}
\usepackage{graphicx}
\usepackage{amsthm}
\usepackage{algorithm}
\usepackage{epic}

\newtheorem{rem}{Remark}
\newenvironment{remark}{\begin{rem}\rm}{\end{rem}}
\newtheorem{exmpl}{Example}
\newenvironment{example}{\begin{exmpl}\rm}{\end{exmpl}}
\newtheorem{theorem}{Theorem}

\def\sign{\mathop{\mathrm{sign}}}
\def\diag{\mathop{\mathrm{diag}}}
\renewcommand{\blacksquare}{\quad \square}

\begin{document}

\begin{frontmatter}

\title{Accurate eigenvalue decomposition of arrowhead matrices and applications}
\author[FESB]{N.\ Jakov\v{c}evi\'{c} Stor\corref{c1}\fnref{f1}}
\ead{nevena@fesb.hr}
\author[FESB]{I.\ Slapni\v{c}ar\fnref{f1}}
\ead{ivan.slapnicar@fesb.hr}
\author[PEN]{J.\ Barlow\fnref{f2}}
\ead{barlow@cse.psu.edu}

\cortext[c1]{Corresponding author}
\fntext[f1]{The research of Ivan Slapni\v{c}ar and Nevena Jakov\v{c}evi\'{c}
  Stor was supported by the Ministry of Science, Education and
  Sports of the Republic of Croatia under grant 023-0372783-1289.} 
\fntext[f2]{The research of Jesse L. Barlow  was supported
by the National Science Foundation under grant  CCF-1115704.}

\address[FESB]{Faculty of Electrical Engineering, Mechanical Engineering and
Naval Architecture, University of Split, Rudjera Bo\v{s}kovi\'{c}a
32, 21000 Split, Croatia}

\address[PEN]{Department of Computer Science and Engineering, The Pennsylvania State
University, University Park, PA 16802-6822, USA}

\begin{abstract}
We present a new algorithm for solving an eigenvalue problem for a real
symmetric arrowhead matrix. The algorithm computes
all eigenvalues and all components of the corresponding eigenvectors with
high relative accuracy in $O(n^{2})$ operations. 
The algorithm is based on a shift-and-invert approach.  
Double precision is eventually needed to compute only one element of
the inverse of the shifted matrix. Each eigenvalue and the corresponding 
eigenvector can be
computed separately, which makes the algorithm adaptable for parallel
computing. Our results extend to Hermitian arrowhead matrices, real
symmetric diagonal-plus-rank-one matrices and singular value decomposition of
real triangular arrowhead matrices.  
\end{abstract}


\end{frontmatter}

\section{Introduction and Preliminaries}

In this paper we consider eigenvalue problem for a real symmetric matrix $A$
which is zero except for its main diagonal and one row and column. Since
eigenvalues are invariant under similarity transformations, we can
symmetrically permute the rows and the columns of the given matrix.
Therefore, we assume without loss of generality that the matrix $A$ is a $%
n\times n$ real symmetric arrowhead matrix of the form

\begin{equation}
A=\left[ 
\begin{array}{cc}
D & z \\ 
z^{T} & \alpha%
\end{array}%
\right],  \label{A}
\end{equation}%
where 
\begin{equation*}
D=\mathop{\mathrm{diag}}(d_{1},d_{2},\ldots ,d_{n-1})\text{ }
\end{equation*}%
is diagonal matrix of order $n-1$,%
\begin{equation}
z=\left[ 
\begin{array}{cccc}
\zeta _{1} & \zeta _{2} & \cdots & \zeta _{n-1}%
\end{array}%
\right] ^{T}  \label{z}
\end{equation}%
is a vector and $\alpha$ is a scalar.

Such matrices arise in the description of radiationless transitions in
isolated molecules \cite{MJ68}, oscillators vibrationally coupled with a
Fermi liquid \cite{Gad81}, quantum optics \cite{MM10} (see also Example 
\ref{ex:quant}). Such matrices also arise in solving symmetric real tridiagonal
eigenvalue problems with the divide-and-conquer method \cite{GE95}.

In this paper we present an algorithm which computes all eigenvalues and all
components of the corresponding eigenvectors with high relative accuracy in $%
O(n^{2})$ operations.

Without loss of generality we may assume that $A$ is irreducible, that is, 
\begin{equation*}
\zeta _{i}\neq 0,\text{ for all }i
\end{equation*}%
and 
\begin{equation*}
d_{i}\neq d_{j},\text{ for all }i\neq j,\text{ }i,j=1,\ldots ,n-1.
\end{equation*}%
If $A$ has a zero in the last column, say $\zeta _{i}=0$, then the diagonal
element $d_{i}$ is an eigenvalue whose corresponding eigenvector is the $i$%
-th unit vector, and we can reduce the size of the problem by deleting the $%
i $-th row and column of the matrix, eventually obtaining a matrix for which
all elements $\zeta _{j}$ are nonzero. If $d_{i}=d_{j}$, then $d_{i}$ is
eigenvalue of matrix $A$ (this follows from the interlacing property (\ref%
{interlace})), and we can reduce the size of the problem by annihilating $%
\zeta_j$ with a Givens rotation in the $(i,j)$-plane and proceeding as in
the previous case.

Further, by symmetric row and column pivoting, we can order elements of $D$
such that 
\begin{equation}
d_{1}>d_{2}>\cdots >d_{n-1}.  \label{order}
\end{equation}%
Hence, we will consider only ordered and irreducible arrowhead matrices.
Without loss of generality we can also assume that $\zeta _{i}>0$ for all $i$%
, which can be attained by pre- and post-multiplication of the matrix $A$
with $D=\mathop{\mathrm{diag}}(\mathop{\mathrm{sign}}(\zeta _{i}))$).

Let 
\begin{equation}
A=V\Lambda V^{T}  \label{Aeigendec}
\end{equation}%
be the eigenvalue decomposition of $A$. Here 
\begin{equation*}
\Lambda =\mathop{\mathrm{diag}}(\lambda _{1},\lambda _{2},\ldots ,\lambda
_{n})
\end{equation*}%
is a diagonal matrix whose diagonal elements are the eigenvalues of $A$, and 
\begin{equation*}
V=\left[ 
\begin{array}{ccc}
v_{1} & \cdots & v_{n}%
\end{array}%
\right]
\end{equation*}%
is an orthonormal matrix whose columns are the corresponding eigenvectors.

The eigenvalues of $A$ are the zeros of the Pick function (see \cite%
{BG92,OlSt90}) 
\begin{equation}
f(\lambda )=\alpha -\lambda -\sum_{i=1}^{n-1}\frac{\zeta _{i}^{2}}{%
d_{i}-\lambda }=\alpha -\lambda -z^{T}(D-\lambda I)^{-1}z,  \label{Pick}
\end{equation}%
and the corresponding eigenvectors are given by 
\begin{equation}
v_{i}=\frac{x_{i}}{\left\Vert x_{i}\right\Vert _{2}},\text{ \ \ \ }x_{i}=%
\left[ 
\begin{array}{c}
\left( D-\lambda _{i}I\right) ^{-1}z \\ 
-1%
\end{array}%
\right] ,\text{ }i=1,\ldots ,n.  \label{Aeigenvec}
\end{equation}%
Diagonal elements of the matrix $D$, $d_{i}$, are called poles of the
function $f$.

Notice that (\ref{order}) and the Cauchy interlacing theorem \cite[Theorem
8.1.7]{GV96} applied to matrices $D$ and $A$ imply the interlacing property 
\begin{equation}
\lambda _{1}>d_{1}>\lambda _{2}>d_{2}>\cdots >d_{n-2}>\lambda
_{n-1}>d_{n-1}>\lambda _{n}.  \label{interlace}
\end{equation}

Since $A$ is symmetric, its eigenvalues may be computed by invoking any of a
number of standard programs (LAPACK \cite{ABB99}). However, these programs
usually begin with an initial reduction of the matrix to tridiagonal form 
\cite{Oli98}, or as proposed in \cite{OlSt90}, with an alternative which
takes advantage of the structure of $A$ by finding the zeros of the Pick
function given in $\left( \ref{Pick}\right) $, for the eigenvalues of $A$.
This results in an algorithm which requires only $O(n^{2})$ computations and 
$O(n)$ storage. Although the idea is conceptually simple and in fact has
been used to solve other eigenvalue problems of special structure 
\cite{Bar93, BN78, Cup81, DS87}, 
the computation is not always stable \cite{GE95}.
Namely, if the computed eigenvalues $\lambda _{i}$ are not accurate enough,
then the computed eigenvectors $v_{i}$ may not be sufficiently orthogonal
(see Example \ref{ex:3}). The existing algorithms for arrowhead matrices 
\cite{GE95,OlSt90} obtain orthogonal eigenvectors with the following
procedure:

\begin{itemize}
\item[-] compute the eigenvalues $\tilde \lambda_i$ of $A$ by solving (\ref%
{Pick});

\item[-] construct a new matrix 
\begin{equation*}
\tilde{A}=\left[ 
\begin{array}{cc}
D & \tilde{z} \\ 
\tilde{z}^{T} & \tilde{\alpha }%
\end{array}%
\right]
\end{equation*}%
by solving inverse problem with the prescribed eigenvalues $\tilde{\lambda }$%
, and diagonal matrix $D$, that is, compute new $\tilde{z}$ and $\tilde{%
\alpha }$ as 
\begin{equation}
\tilde{\zeta}_{i} =\sqrt{\left( d_{i}-\tilde{\lambda }_{n}\right) \left( 
\tilde{\lambda }_{1}-d_{i}\right) \prod\limits_{j=2}^{i}\frac{\left( \tilde{%
\lambda }_{j}-d_{i}\right) }{\left( d_{j-1}-d_{i}\right) }%
\prod\limits_{j=i+1}^{n-1}\frac{\left( \tilde{\lambda }_{j}-d_{i}\right) }{%
\left( d_{j}-d_{i}\right) }}\text{,}  \notag
\end{equation}%
\begin{equation}
\tilde{\alpha }=\tilde{\lambda }_{n}+\sum\limits_{j=1}^{n-1}\left( \tilde{%
\lambda }_{j}-d_{j}\right) \text{.}  \notag
\end{equation}

\item[-] compute eigenvectors of $\tilde{A}$ by (\ref{Aeigenvec}).
\end{itemize}

Since the formulas for $\tilde{\zeta_i}$ involve only multiplications,
division and subtractions of exact quantities, each $\zeta_i$ is computed
with relative error of $O(\varepsilon_M)$, where $\varepsilon_M$ denotes the
machine precision.\footnote{%
The machine precision $\varepsilon_M$ is defined as a smallest positive
number such that in the floating-point arithmetic $1+\varepsilon_M\ne 1$. In
Matlab or FORTRAN REAL(8) arithmetic $\varepsilon_M= 2.2204\cdot 10^{-16}$,
thus the floating-point numbers have approximately 16 significant decimal
digits. The term ``double of the working precision'' means that the
computations are performed with numbers having approximately 32 significant
decimal digits, or with the machine precision equal to $\varepsilon_M^2$.}
Therefore, $\tilde A=A+\delta A$, where $\| \delta A\|_2=O(\epsilon_M)$.
Here $\| \cdot \|_2$ denotes the spectral matrix norm. We conclude that the
computed eigenvalues $\tilde \lambda_i$ satisfy standard perturbation bounds
like those from \cite[Corollary 8.1.6]{GV96}. Further, since $\tilde
\lambda_i$ are the eigenvalues of the matrix $\tilde A$ computed to higher
relative accuracy, the eigenvectors computed by (\ref{Aeigenvec}) are
orthogonal to machine precision. For details see \cite{GE95,OlSt90}.

Our algorithm uses a different approach. Accuracy of the eigenvectors and
their orthogonality follows from high relative accuracy of the computed
eigenvalues and there is no need for follow-up orthogonalization. 
The algorithm is based on shift-and-invert technique. Basically, the
eigenvalue $\lambda $ is computed as the largest or the smallest eigenvalue
of the inverse of the matrix shifted to the pole $d_{i}$ which is nearest to 
$\lambda $, that is, 
\begin{equation}  \label{8}
\lambda=\frac{1}{\nu}+d_i,
\end{equation}
where $\nu$ is either smallest or largest eigenvalue of the matrix 
\begin{equation*}
A_i^{-1}\equiv (A-d_iI)^{-1}.
\end{equation*}

Inverses of arrowhead matrices are structured in the following manner (here $%
\times$ stands for non-zero element): the inverse of an arrowhead matrix
with zero on the shaft is a permuted arrowhead matrix with zero on the
shaft, 
\begin{align*}
\begin{bmatrix}
\times &  &  &  & \times \\ 
& \times &  &  & \times \\ 
&  & 0 &  & \times \\ 
&  &  & \times & \times \\ 
\times & \times & \times & \times & \times%
\end{bmatrix}%
^{-1} &= 
\begin{bmatrix}
\times &  & \times &  &  \\ 
& \times & \times &  &  \\ 
\times & \times & \times & \times & \times \\ 
&  & \times & \times &  \\ 
&  & \times &  & 0%
\end{bmatrix}%
,
\end{align*}
and the inverse of the full arrowhead matrix is a diagonal-plus-rank-one
(DPR1) matrix, 
\begin{align*}
\begin{bmatrix}
\times &  &  &  & \times \\ 
& \times &  &  & \times \\ 
&  & \times &  & \times \\ 
&  &  & \times & \times \\ 
\times & \times & \times & \times & \times%
\end{bmatrix}%
^{-1} &= 
\begin{bmatrix}
\times &  &  &  &  \\ 
& \times &  &  &  \\ 
&  & \times &  &  \\ 
&  &  & \times &  \\ 
&  &  &  & 0%
\end{bmatrix}
\pm u u^T.
\end{align*}

Our algorithm is completely parallel, since the computation of one
eigenvalue and its eigenvector is completely independent of the computation
of other eigenvalues and eigenvectors.

In Section 2 we describe the basic idea of our algorithm named $aheig$
(ArrowHead EIGenvalues). In Section 3 we discuss the accuracy of the
algorithm. In Section 4 we present the complete algorithm which uses double
of the working precision, if necessary. In Section 5 we illustrate algorithm
with few examples and in Section 6 we apply our results to 
eigenvalue decomposition of Hermitian arrowhead matrix, singular value
decomposition of real triangular arrowhead matrix and  
eigenvalue decomposition of real symmetric diagonal-plus-rank-one matrix. 
The proofs are given in Appendix A.

\section{Basic shift-and-invert algorithm}

Let $\lambda $ be an eigenvalue of $A$, let $v$ be its eigenvector, and let $%
x$ be the unnormalized version of $v$ from (\ref{Aeigenvec}). Let $d_{i}$ be
the pole which is closest to $\lambda $. Clearly, from (\ref{interlace}) it
follows that either $\lambda =\lambda _{i} $ or $\lambda =\lambda _{i+1}$.
Let $A_{i}$ be the shifted matrix 
\begin{equation}
A_{i}=A-d_{i}I=\left[ 
\begin{array}{cccc}
D_{1} & 0 & 0 & z_{1} \\ 
0 & 0 & 0 & \zeta _{i} \\ 
0 & 0 & D_{2} & z_{2} \\ 
z_{1}^{T} & \zeta _{i} & z_{2}^{T} & a%
\end{array}%
\right] \text{,}  \label{Ai}
\end{equation}%
where%
\begin{align*}
D_{1}& =\mathop{\mathrm{diag}}(d_{1}-d_{i},\ldots ,d_{i-1}-d_{i}), \\
D_{2}& =\mathop{\mathrm{diag}}(d_{i+1}-d_{i},\ldots ,d_{n-1}-d_{i}), \\
z_{1}& =\left[ 
\begin{array}{cccc}
\zeta _{1} & \zeta _{2} & \cdots & \zeta _{i-1}%
\end{array}%
\right] ^{T}, \\
z_{2}& =\left[ 
\begin{array}{cccc}
\zeta _{i+1} & \zeta _{i+2} & \cdots & \zeta _{n-1}%
\end{array}%
\right] ^{T}, \\
a& =\alpha -d_{i}.
\end{align*}%
Notice that $D_{1}$ $(D_{2})$ is positive (negative) definite.

Obviously, if $\lambda $ is an eigenvalue of $A$, then 
\begin{equation*}
\mu =\lambda -d_{i}
\end{equation*}%
is an eigenvalue of $A_{i}$, and vice versa, and they both have the same
eigenvector. The inverse of $A_{i}$ is%
\begin{equation}
A_{i}^{-1}=\left[ 
\begin{array}{cccc}
D_{1}^{-1} & w_{1} & 0 & 0 \\ 
w_{1}^{T} & b & w_{2}^{T} & 1/\zeta _{i} \\ 
0 & w_{2} & D_{2}^{-1} & 0 \\ 
0 & 1/\zeta _{i} & 0 & 0%
\end{array}%
\right] ,  \label{invAi}
\end{equation}%
where%
\begin{align}
w_{1}& =-D_{1}^{-1}z_{1}\frac{1}{\zeta _{i}},  \notag \\
w_{2}& =-D_{2}^{-1}z_{2}\frac{1}{\zeta _{i}},  \notag \\
b& =\frac{1}{\zeta _{i}^{2}}\left(
-a+z_{1}^{T}D_{1}^{-1}z_{1}+z_{2}^{T}D_{2}^{-1}z_{2}\right) .  \label{b_ob}
\end{align}

Notice that 
\begin{equation*}
b=\bar{f}\left( d_{i}\right) /\zeta _{i}^{2}\text{ }
\end{equation*}%
where%
\begin{equation*}
\bar{f}\left( d_{i}\right) =\alpha -d_{i}-\bar{z}^{T}\left( \bar{D}%
-d_{i}I\right) ^{-1}\bar{z}
\end{equation*}%
where $\bar{D}$ is the diagonal matrix $D$ without $d_{i}$ and $\bar{z}$ is $%
z$ without $\zeta _{i}$.

The eigenvector $x$ from (\ref{Aeigenvec}) is given by 
\begin{equation}
x=\begin{bmatrix} x_1\\ \vdots \\ x_n
\end{bmatrix}
=
\left[ 
\begin{array}{c}
\left( D_{1}-\mu I\right) ^{-1}z_{1} \\ 
-\displaystyle\frac{\zeta _{i}}{\mu } \\ 
\left( D_{2}-\mu I\right) ^{-1}z_{2} \\ 
-1%
\end{array}%
\right] .  \label{eigenvecAi}
\end{equation}

If $\lambda $ is an eigenvalue of $A$ which is closest to the pole $d_{i}$,
then $\mu $ is the eigenvalue of matrix $A_{i}$ which is closest to zero and 
\begin{equation*}
\nu=\frac{1}{\mu}=\pm \left\Vert A_{i}^{-1}\right\Vert _{2}.
\end{equation*}
In this case, if all entries of $A_{i}^{-1}$ are computed with high relative
accuracy, then, according to standard perturbation theory, $\nu$ is computed
to high relative accuracy (by any reasonable algorithm). In Section 3 we
show that all entries of $A_{i}^{-1}$ are indeed computed to high relative
accuracy, except possibly $b$ (see (\ref{b_ob})). If $b$ is not computed to
high relative accuracy and it influences $\left\Vert A_{i}^{-1}\right\Vert
_{2}$, it is sufficient to compute it in double of the working precision
(see Section 4).


Further, if $\mu $ is not the eigenvalue of $A_{i}$ which is closest to
zero, then $\left\vert \nu \right\vert <\left\Vert A_{i}^{-1}\right\Vert
_{2} $, and the quantity 
\begin{equation}
K_{\nu} = \frac{\left\Vert A_{i}^{-1}\right\Vert _{2}}{\left\vert \nu
\right\vert}  \label{K0}
\end{equation}%
tells us how far is $\nu $ from the absolutely largest eigenvalue of $%
A_{i}^{-1}$. If $K_{\nu}\gg 1$, then the standard perturbation theory does
not guarantee that the eigenvalue $\mu $ will be computed with high relative
accuracy. Remedies of this situation are described in Remark \ref{rem:2}.


With this approach the componentwise high relative accuracy of the
eigenvectors computed by $\left( \ref{eigenvecAi}\right) $ follows from high
relative accuracy of the computed eigenvalues (see Theorem \ref{T2}).
Componentwise high relative accuracy of the computed eigenvectors implies,
in turn, their orthogonality.

The described procedure is implemented in algorithm $aheig\_basic$
(Algorithm 1). The computation of the inverse of the shifted matrix, $%
A_{i}^{-1}$, according to formulas (\ref{invAi}) and (\ref{b_ob}), is
implemented in Algorithm 2. Algorithm 3 computes the largest or the smallest
zero of the Pick function (\ref{Pick}) by bisection. Given eigenvalue $%
\lambda$, Algorithm 4 computes the corresponding eigenvector by (\ref%
{Aeigenvec}) or (\ref{eigenvecAi}), respectively.

\begin{algorithm}
\caption{} $[\lambda,v]=\mathbf{aheig\_basic}\left(
D,z,\alpha ,k\right)$

{\small \% Computes the} $k${\small -th eigenpair
of an irreducible arrowhead matrix }

{\small \%} $A=[\diag\left( D\right) \ z;
z^{\prime } \ \alpha ]$




$n=max(size(D))+1$



{\small \% Determine the shift} $\sigma${\small , the shift index}
$i${\small , and whether} $\lambda$ {\small is
  on the left}

{\small \%  or the right side of the nearest pole.}

{\small \% Exterior eigenvalues (}$k=1${\small \ or  }$k=n${\small ):}%

\textbf{if }$k==1$









\ \ $\sigma=d_{1}$

\ \ $i=1$

\ \ $side=\mathrm{'R'}$


\textbf{elseif } $k==n$








\ \ $\sigma=d_{n-1}$

\ \ $i=n-1$

\ \ $side=\mathrm{'L'}$





\textbf{else }

\ \ {\small \% Interior eigenvalues (}$k\in\{2,\ldots,n-1\}${\small ):}


\ \ $Dtemp=D-d_{k}$

\ \ $atemp=\alpha -d_{k}$

\ \ $middle=Dtemp_{k-1}/2$

\ \ $Fmiddle=atemp-middle-\sum (z^{2}./(Dtemp-middle))$

\ \ \textbf{if }$Fmiddle<0$

\ \ \ \ $\sigma=d_{k}$

\ \ \ \ $i=k$

\ \ \ \ $side=\mathrm{'R'}$

\ \ \textbf{else}

\ \ \ \ $\sigma=d_{k-1}$

\ \ \ \ $i=k-1$

\ \ \ \ $side=\mathrm{'L'}$

\ \ \textbf{end}

\textbf{end}


{\small \% Compute the inverse of the shifted matrix, }$A_{i}^{-1}$

$[invD_1, invD_2,w_1, w_2, w_\zeta, b]=\mathbf{invA}(D,z,\alpha ,i)$

{\small \% Compute the leftmost or the rightmost eigenvalue of }$A_{i}^{-1}$%

$\nu=\mathbf{bisect}([invD_1;0;invD_2],[w_1; w_\zeta;w_2],b,side)$

{\small \% Compute the corresponding eigenvector}

$\mu=1/\nu$

$v=\mathbf{vect}(D-\sigma,z,\mu)$

{\small \% Shift the eigenvalue back}

$\lambda=\mu+\sigma$


\end{algorithm}

\begin{algorithm}

\caption{} $[invD_1, invD_2,w_1, w_2, w_\zeta, b ]=\mathbf{invA}\left( D,z,\alpha,i\right)$

{\small \% Computes the inverse of an arrowhead matrix}
$A=[\diag(D-d_{i}) \ z; z^{\prime } \ \alpha -d_{i}]$ 

{\small \% according to (\ref{invAi}) and (\ref{b_ob}).}


$n=max(size(D))+1$

$D=D-d_{i}$

$a=\alpha -d_{i}$

$w_{1}=-z_{1:i-1}./D_{1:i-1}/z_{i}$

$w_{2}=-z_{i+1:n-1}./D_{i+1:n-1}/z_{i}$


$w_\zeta=1/z_{i}$

$invD_1=1./D_{1:i-1}$

$invD_2=1./D_{i+1:n-1}$

$b =(-a+ sum(z_{1:i-1}.\verb|^| 2./D_{1:i-1}) + sum
(z_{i+1:n-1}.\verb|^|2./D_{i+1:n-1})) /z_{i}\verb|^|2$

\end{algorithm}

\begin{algorithm}
\caption{} $\lambda=\mathbf{bisect}\left( D,z,\alpha,side\right)$

{\small \% Computes the leftmost (for $side$='L') or the rightmost (for
  $side$='R') eigenvalue}

{\small \% } of
an arrowhead matrix $A=[\diag\left( D\right) ${\small \ }$z;${\small \ }$z^{\prime
} ${\small \ }$\alpha ]${\small \ by bisection.}



$n=max(size(D))+1$

{\small \% Determine the starting interval for bisection,} $[left,right]$

\textbf{if} $side==\mathrm{'L'}$


\ \ $left=\min\{D-|z|,\alpha-\|z\|_1\}$

\ \ $right=\min d_{i}$


\textbf{else}

\ \ $right=\max\{D+|z|,\alpha+\|z\|_1\}$

\ \ $left=\max d_{i}$


\textbf{end}

{\small \%\ Bisection}


$middle=(left+right)/2$

\textbf{while}\ $(right-left)/abs(middle)>2\ast eps$

\ \ $Fmiddle=\alpha-middle-sum(z.\verb|^|2./\ (D-middle))$

\ \ \textbf{if\ }\ $Fmiddle>0$

\ \ \ \ $left=middle$

\ \ $\textbf{else}$

\ \ \ \ $right=middle$

\ \ \textbf{end}

\ \ $middle=(left+right)/2$


\textbf{end}

{\small \%\ Eigenvalue }\textit{\ }

$\lambda =right$

\end{algorithm}

\begin{algorithm}
\caption{} $v=\mathbf{vect}\left( D,z,\lambda\right)$

{\small \% Computes the eigenvector of an arrowhead matrix}
$A=[\diag( D) \ z; z^{\prime} \ \alpha]$

{\small \% which corresponds to the eigenvalue} $\lambda$
{\small by using (\ref{Aeigenvec}).}



$v=[ z ./ (D-\lambda);-1]$

$v=v/\|v\|_2$

\end{algorithm}

\section{Accuracy of the algorithm}

We now consider numerical properties of Algorithms 1, 2, 3, and 4. 
We assume tha standard model of floating point arithmetic where subtraction
is preformed with guard digit, such that \cite{Gol91,Par80,GV96,Wil65}
$$
fl (a\circ b) = (a\circ b)(1+\varepsilon_\circ), \quad  
|\varepsilon_\circ |\leq \varepsilon_M, \quad \circ \in\{+,-,*,/\},
$$
where $\varepsilon_M$ is machine precision. In the statements of the theorems
and their proofs we shall use the standard first order approximations, that is,
we neglect the terms of order $O(\varepsilon_M^2)$ or higher. Moreover, we
assume that neither overflow or underflow occurs during the computation.

We shall use
the following notation: 
\begin{equation}\label{tablica}
\begin{tabular}{c|c|c}
Matrix & Exact eigenvalue & Computed eigenvalue \\ \hline
$A$ & $\lambda$ & $\widetilde{\lambda}$ \\ \hline
$A_{i}$ & $\mu $ & $-$ \\ \hline
$\widetilde{A}_{i}=fl(A_{i})$ & $\widehat{\mu }$ & $\widetilde{\mu }=fl(%
\widehat{\mu })$ \\ \hline
$A_{i}^{-1}$ & $\nu $ & $-$ \\ \hline
$\widetilde{(A_{i}^{-1})}=fl(A_{i}^{-1})$ & $\widehat{\nu }$ & $\widetilde{%
\nu }=fl(\widehat{\nu })$ \\ \hline
\end{tabular}%
\end{equation}%
%
%
%
%
%
Here 
\begin{equation*}
\widetilde{A}_{i}=fl\left( A_{i}\right) = \left[ 
\begin{array}{cccc}
D_{1}\left( I+E_{1}\right) & 0 & 0 & z_{1} \\ 
0 & 0 & 0 & \zeta _{i} \\ 
0 & 0 & D_{2}\left( I+E_{2}\right) & z_{2} \\ 
z_{1}^{T} & \zeta _{i} & z_{2}^{T} & a\left( 1+\varepsilon_{a}\right)%
\end{array}%
\right] ,  
\end{equation*}%
where $E_{1}$ and $E_{2}$ are diagonal matrices whose elements are bounded
by $\varepsilon_{M}$ in absolute values and $\left\vert \varepsilon
_{a}\right\vert \leq \varepsilon_{M}$.

Further we define the quantities $\kappa _{\lambda }$, $\kappa _{\mu }$ and $%
\kappa _{b}$ 
as follows: 
\begin{align}
\widetilde{\lambda }&=fl\left( \lambda\right) =\lambda\left( 1+\kappa
_{\lambda }\varepsilon _{M}\right) ,  \label{kapa_lambda} \\
\widetilde{\mu }&=fl\left( \mu \right) =\mu \left( 1+\kappa _{\mu
}\varepsilon _{M}\right) ,  \label{kapa_mi} \\
\widetilde{b}&=fl\left( b\right) =b\left( 1+\kappa _{b}\varepsilon
_{M}\right).  \label{kapa_b}
\end{align}
We also define the quantity 
\begin{equation}
K_{b}=\frac{\left\vert a\right\vert +\left\vert
z_{1}^{T}D_{1}^{-1}z_{1}\right\vert +\left\vert
z_{2}^{T}D_{2}^{-1}z_{2}\right\vert }{\left\vert
-a+z_{1}^{T}D_{1}^{-1}z_{1}+z_{2}^{T}D_{2}^{-1}z_{2}\right\vert }.
\label{K2}
\end{equation}

\subsection{Connection between accuracy of $\protect\lambda$ and $\protect%
\mu $}

Let 
\begin{equation*}
\lambda=\mu +d_{i}
\end{equation*}%
be an eigenvalue of the matrix $A$, where $\mu$ is the corresponding
eigenvalue of the shifted matrix $A_{i}=A-d_{i}$ from which $\lambda$ is
computed. Let 
\begin{equation*}
\widetilde\lambda=fl(\widetilde\mu +d_{i})
\end{equation*}%
be the computed eigenvalue. Theorem \ref{T1} gives us dependency of accuracy
of $\widetilde{\lambda}$ in (\ref{kapa_lambda}) upon accuracy of $\widetilde{%
\mu }$ in (\ref{kapa_mi}).

\begin{theorem}
\label{T1} For $\lambda$ and $\widetilde{\lambda}$ from (\ref{kapa_lambda})
and $\mu$ and $\widetilde{\mu }$ from (\ref{kapa_mi}) we have 
\begin{equation}
\left\vert \kappa _{\lambda }\right\vert \leq \frac{\left\vert
d_{i}\right\vert +\left\vert \mu \right\vert }{\left\vert \lambda
\right\vert }\left( \left\vert \kappa _{\mu }\right\vert +1\right).
\label{lambda_tilda}
\end{equation}
\end{theorem}

Proofs of this theorem and subsequent theorems are given in
Appendix A.

From Theorem \ref{T1} we see that the accuracy of $\widetilde\lambda$
depends on $\kappa _{\mu }$ and the size of the quotient 
\begin{equation}
\frac{\left\vert d_{i}\right\vert +\left\vert \mu \right\vert }{\left\vert
\lambda \right\vert }\text{.}  \label{eq:quotient}
\end{equation}%
%
%
%
%
%
Theorem \ref{C1} analyzes the quotient (\ref{eq:quotient}) with respect to
the position of $\lambda $ and signs of $\mu $ and the neighboring poles.

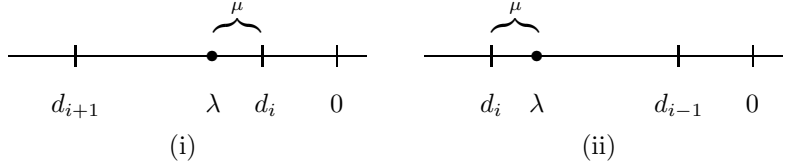
\begin{figure}[hbtp]
\begin{tabular}{cc}
\setlength{\unitlength}{0.00087489in}
\begingroup\makeatletter\ifx\SetFigFont\undefined%
\gdef\SetFigFont#1#2#3#4#5{%
  \reset@font\fontsize{#1}{#2pt}%
  \fontfamily{#3}\fontseries{#4}\fontshape{#5}%
  \selectfont}%
\fi\endgroup%
{\renewcommand{\dashlinestretch}{30}
\begin{picture}(2154,852)(0,-10)
\put(1227,390){\circle*{60}}
\drawline(402,455)(402,325)
\drawline(1977,455)(1977,325)
\drawline(2,390)(2152,390)
\drawline(1527,455)(1527,325)
\put(1482,50){\makebox(0,0)[lb]{\smash{{\SetFigFont{10}{16.8}{\rmdefault}
{\mddefault}{\updefault}$d_i$}}}}
\put(1935,50){\makebox(0,0)[lb]{\smash{{\SetFigFont{10}{16.8}{\rmdefault}
{\mddefault}{\updefault}$0$}}}}
\put(257,50){\makebox(0,0)[lb]{\smash{{\SetFigFont{10}{16.8}{\rmdefault}
{\mddefault}{\updefault}$d_{i+1}$}}}}
\put(1187,50){\makebox(0,0)[lb]{\smash{{\SetFigFont{10}{16.8}{\rmdefault}
{\mddefault}{\updefault}$\lambda$}}}}
\put(1215,550){\makebox(0,0)[lb]{\smash{{\SetFigFont{10}{16.8}{\rmdefault}
{\mddefault}{\updefault}$\mu
\atop \overbrace{\phantom{x}}$}}}}
\end{picture}
}
& 
\setlength{\unitlength}{0.00087489in}
\begingroup\makeatletter\ifx\SetFigFont\undefined%
\gdef\SetFigFont#1#2#3#4#5{%
  \reset@font\fontsize{#1}{#2pt}%
  \fontfamily{#3}\fontseries{#4}\fontshape{#5}%
  \selectfont}%
\fi\endgroup%
{\renewcommand{\dashlinestretch}{30}
\begin{picture}(2154,852)(0,-10)
\put(675,390){\circle*{60}}
\drawline(402,455)(402,325)
\drawline(1977,455)(1977,325)
\drawline(2,390)(2152,390)
\drawline(1527,455)(1527,325)
\put(1382,50){\makebox(0,0)[lb]{\smash{{\SetFigFont{10}{16.8}{\rmdefault}
{\mddefault}{\updefault}$d_{i-1}$}}}}
\put(1935,50){\makebox(0,0)[lb]{\smash{{\SetFigFont{10}{16.8}{\rmdefault}
{\mddefault}{\updefault}$0$}}}}
\put(335,50){\makebox(0,0)[lb]{\smash{{\SetFigFont{10}{16.8}{\rmdefault}
{\mddefault}{\updefault}$d_{i}$}}}}
\put(630,50){\makebox(0,0)[lb]{\smash{{\SetFigFont{10}{16.8}{\rmdefault}
{\mddefault}{\updefault}$\lambda$}}}}
\put(385,550){\makebox(0,0)[lb]{\smash{{\SetFigFont{10}{16.8}{\rmdefault}
{\mddefault}{\updefault}$\mu
\atop \overbrace{\phantom{x}}$}}}}
\end{picture}
}
\\ 
(i) & (ii)%
\end{tabular}
\caption{Typical situations from Theorem \protect\ref{C1}}
\label{K1sl}
\end{figure}

\begin{theorem}
\label{C1} Let the assumptions of Theorem \ref{T1} hold.

\begin{itemize}
\item[(i)] If (see Figure \ref{K1sl} (i)) 
\begin{equation*}
\mathop{\mathrm{sign}}\left( d_{i}\right) =\mathop{\mathrm{sign}}\left( \mu
\right),
\end{equation*}
then 
\begin{equation*}
\frac{\left\vert d_{i}\right\vert +\left\vert \mu \right\vert } {\left\vert
\lambda \right\vert }= 1.
\end{equation*}

\item[(ii)] If $\lambda $ is between two poles of the same sign and $%
\mathop{\mathrm{sign}}\left( d_{i}\right) \neq \mathop{\mathrm{sign}}\left(
\mu \right) $ (see Figure \ref{K1sl} (ii)), then 
\begin{equation*}
\frac{\left\vert d_{i}\right\vert +\left\vert \mu \right\vert }{\left\vert
\lambda \right\vert }\leq 3.
\end{equation*}
\end{itemize}
\end{theorem}

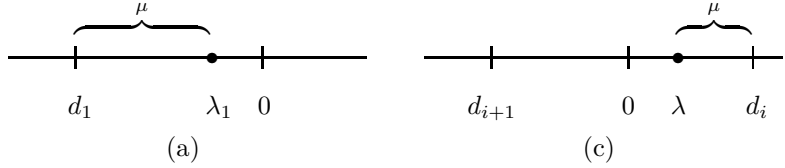
\begin{figure}[hbtp]
\begin{tabular}{cc}
\setlength{\unitlength}{0.00087489in} \begingroup\makeatletter\ifx\SetFigFont%
\undefined\gdef\SetFigFont#1#2#3#4#5{\ \reset@font\fontsize{#1}{#2pt} %
\fontfamily{#3}\fontseries{#4}\fontshape{#5} \selectfont}\fi%
\endgroup{\renewcommand{\dashlinestretch}{30}
\begin{picture}(2154,852)(0,-10)
\put(1227,390){\circle*{60}}
\drawline(402,455)(402,325)
\drawline(2,390)(2152,390)
\drawline(1527,455)(1527,325)
\put(1500,50){\makebox(0,0)[lb]{\smash{{\SetFigFont{10}{16.8}{\rmdefault}
{\mddefault}{\updefault}$0$}}}}
\put(357,50){\makebox(0,0)[lb]{\smash{{\SetFigFont{10}{16.8}{\rmdefault}
{\mddefault}{\updefault}$d_{1}$}}}} \put(1187,50){\makebox(0,0)[lb]{\smash{{\SetFigFont{10}{16.8}{\rmdefault}
{\mddefault}{\updefault}$\lambda_1$}}}} \put(375,550){\makebox(0,0)[lb]{\smash{{\SetFigFont{10}{16.8}{\rmdefault}
{\mddefault}{\updefault}$\mu
\atop \overbrace{\phantom{xxxxxxxxx}}$}}}} \end{picture} } & %
\setlength{\unitlength}{0.00087489in} \begingroup\makeatletter\ifx\SetFigFont%
\undefined\gdef\SetFigFont#1#2#3#4#5{\ \reset@font\fontsize{#1}{#2pt} %
\fontfamily{#3}\fontseries{#4}\fontshape{#5} \selectfont}\fi%
\endgroup{\renewcommand{\dashlinestretch}{30}
\begin{picture}(2154,652)(0,-10)
\put(1527,390){\circle*{60}}
\drawline(402,455)(402,325)
\drawline(1977,455)(1977,325)
\drawline(2,390)(2152,390)
\drawline(1227,455)(1227,325)
\put(1482,50){\makebox(0,0)[lb]{\smash{{\SetFigFont{10}{16.8}{\rmdefault}
{\mddefault}{\updefault}$\lambda$}}}}
\put(1935,50){\makebox(0,0)[lb]{\smash{{\SetFigFont{10}{16.8}{\rmdefault}
{\mddefault}{\updefault}$d_i$}}}} \put(257,50){\makebox(0,0)[lb]{\smash{{\SetFigFont{10}{16.8}{\rmdefault}
{\mddefault}{\updefault}$d_{i+1}$}}}} \put(1187,50){\makebox(0,0)[lb]{\smash{{\SetFigFont{10}{16.8}{\rmdefault}
{\mddefault}{\updefault}$0$}}}} \put(1495,550){\makebox(0,0)[lb]{\smash{{\SetFigFont{10}{16.8}{\rmdefault}
{\mddefault}{\updefault}$\mu
\atop \overbrace{\phantom{xxxxx}}$}}}} \end{picture} } \\ 
(a) & (c)%
\end{tabular}
\caption{Typical situations for special cases}
\label{K1_ac_sl}
\end{figure}

Theorem \ref{C1} does not cover the following cases:

\begin{itemize}
\item[(a)] If $d_{1}<0$, then $\mu >0$. If, further, $\left\vert
d_{1}\right\vert \approx \left\vert \mu \right\vert$, then $\lambda _{1}$ is
near zero, and $(\left\vert d_{1}\right\vert +\left\vert \mu
\right\vert)/\left\vert \lambda _{1}\right\vert\gg 1$ (see Figure \ref%
{K1_ac_sl} (a)).\footnote{%
In this case $\lambda _{1}$ is computed as a difference of two close
quantities and cancellation can occur.}

\item[(b)] If $d_{n}>0,$ then $\mu <0$. If, further, $\left\vert
d_{n}\right\vert \approx \left\vert \mu \right\vert$, then $\lambda _{n}$ is
near zero, and again $(\left\vert d_{n}\right\vert +\left\vert \mu
\right\vert)/\left\vert \lambda _{n}\right\vert\gg 1$.\footnote{%
In this case $\lambda _{n}$ is computed as a difference of two close
quantities and cancellation can occur.}

\item[(c)] If $\lambda $ is between two poles of the different signs
and $\mathop{\mathrm{sign}}\left( d_{i}\right) \neq \mathop{\mathrm{sign}}%
\left( \mu \right) $, then either $d_{i+1}<0<d_{i}$ and $\mu <0$, or $%
d_{i}<0<d_{i-1}$ and $\mu >0$. In both cases, if, additionally, $\left\vert
d_{i}\right\vert \approx \left\vert \mu \right\vert $, then $\lambda $ is
near zero, and $(\left\vert d_{i}\right\vert +\left\vert \mu
\right\vert)/\left\vert \lambda \right\vert\gg 1$ (see Figure \ref{K1_ac_sl}
(c)).
\end{itemize}

Since only one of these three cases can occur, Theorems \ref{T1} and \ref{C1}
imply 
that for all eigenvalues $\lambda \in \sigma \left( A\right)$, but
eventually one, it holds 
\begin{equation*}
\frac{\left\vert d_{i}\right\vert +\left\vert \mu \right\vert }{\left\vert
\lambda \right\vert }\leq 3.  
\end{equation*}

If one of the above cases does occur, remedies are given in the following
remark.

\begin{remark}
\label{rem:1} If one of the cases (a), (b) or (c) 
occurs, then $\lambda$ is an
eigenvalue of $A$ nearest to zero, and we can accurately compute it from the
inverse of $A$. Notice that the inverse is of an unreduced arrowhead matrix
with non-zero shaft is a diagonal-plus-rank-one (DPR1) matrix of the form 
\begin{equation*}
A^{-1}=%
\begin{bmatrix}
D^{-1} &  \\ 
& 0%
\end{bmatrix}%
+\rho uu^{T},
\end{equation*}%
where 
\begin{equation*}
u=%
\begin{bmatrix}
z^{T}D^{-1} & -1%
\end{bmatrix}%
^{T},\quad \rho =\frac{1}{a-z^{T}D^{-1}z}.
\end{equation*}%
Eigenvalues of $A^{-1}$ are zeros of (see \cite{Bar93,Mel95}) 
\begin{equation*}
\varphi \left( \lambda \right) =1+\rho \sum_{j=1}^{n}\frac{u_{j}^{2}}{%
d_{j}-\lambda }.
\end{equation*}%
Since the absolutely largest eigenvalue of $A^{-1}$ is computed accurately
according to standard perturbation theory, and $1/|\lambda|=\|A^{-1}\|_2$, $%
\lambda$ is also computed with high relative accuracy. In computing matrix $%
A^{-1}$, eventually $\rho $ needs to be computed in higher precision. For
more details see Remark \ref{rem:2}.
If the denominator in $\rho$ is computed as zero, the matrix $A$ is
numerically singular and we can set $\lambda=0$. Notice that all components of
the the corresponding eigenvector are still computed accurately. 
\end{remark}

\begin{remark}\label{rem:1a}
Notice that Algorithm 1 (and, consequently, Algorithm 5 below) 
can be easily modified to return both quantities,
$d_i$ and $\mu$ such that $\lambda=d_i+\mu$. If none of the remedies from
Remark \ref{rem:1} were needed, these two quantities give additional information
about $\lambda$ (that is, they give a more accurate representation of
$\lambda$). An example is given in Example \ref{ex:2}.
\end{remark}

We still need to bound the quantity $\kappa _{\mu }$ from (\ref{lambda_tilda}%
). This quantity essentially depends on the accuracy of $fl(b)$. 
The bound for $\kappa_{\mu }$ is given in Theorem \ref{T7}.

\subsection{Accuracy of the eigenvectors}

Since the eigenvector is computed by (\ref{eigenvecAi}), its accuracy
depends on the accuracy of $\widetilde\mu$ as described by the following
theorem:

\begin{theorem}
\label{T2} Let (\ref{kapa_mi}) hold and let 
\begin{equation}
\widetilde{x}=\left[ 
\begin{array}{c}
\widetilde{x}_{1} \\ 
\vdots \\ 
\\ 
\widetilde{x}_{n}%
\end{array}%
\right] =fl(\left[ 
\begin{array}{c}
\left( D_{1}\left( I+E_{1}\right) -\widetilde{\mu }I\right) ^{-1}z_{1} \\ 
-\displaystyle\frac{\zeta _{i}}{\widetilde{\mu }} \\ 
\left( D_{2}\left( I+E_{2}\right) -\widetilde{\mu }I\right) ^{-1}z_{2} \\ 
-1%
\end{array}%
\right] )  \label{x_tilda}
\end{equation}%
be the computed un-normalized eigenvector corresponding to $\mu$ and $%
\lambda $. Then 
\begin{equation*}
\widetilde{x_{j}}=x_{j}\left( 1+\varepsilon _{x_{j}}\right) \text{, \ }%
\left\vert \varepsilon _{x_{j}}\right\vert \leq 3\left( |\kappa _{\mu
}|+3\right) \varepsilon _{M}, \quad j=1,\ldots ,n.
\end{equation*}
\end{theorem}

In other words, if $\kappa_\mu$ is small, then all components of the
eigenvector are computed to high relative accuracy. Since the accuracy of $%
\widetilde{\lambda }$ and $\widetilde{x}$ depends on the accuracy of $%
\widetilde{\mu }$ (on the size of $\kappa_\mu$) in the next three
subsections tells we discuss the accuracy of $\widetilde{\mu}$. Since $%
\widetilde{\mu}$ is computed as an inverse of the eigenvalue of the matrix $%
fl(A_{i}^{-1})$, we first discuss the accuracy of that matrix.

\subsection{Accuracy of the matrix $A_{i}^{-1}$}

We have the following theorem:

\begin{theorem}
\label{T3} For the computed elements of the matrix $A_{i}^{-1}$ from (\ref%
{invAi}) and (\ref{b_ob}) for all $(j,k)\neq (i,i)$ we have 
\begin{equation*}
\widetilde{\left( A_{i}^{-1}\right) }_{jk}=fl\left( A_{i}^{-1}\right)
_{jk}=\left( A_{i}^{-1}\right) _{jk}(1+\varepsilon _{jk}),\text{ }\left\vert
\varepsilon _{jk}\right\vert \leq 3\varepsilon _{M}.
\end{equation*}%
For the computed element $b\equiv\left( A_{i}^{-1}\right) _{ii}$ from (\ref%
{kapa_b}) we have 
\begin{equation*}
\left\vert \kappa _{b}\right\vert \leq (n+3)K_{b},  
\end{equation*}
where $K_b$ is defined by (\ref{K2}).
\end{theorem}


The above theorem states that all elements of the matrix $A_{i}^{-1}$ are
computed with high relative accuracy except possibly $b$. Therefore, we have
to monitor whether $b$ is computed accurately, and, if not, it needs to be
computed in double of the working precision (see Section 4 for details).

\subsection{Accuracy of bisection}

Let $\lambda _{\max }$ be the absolutely largest eigenvalue of a symmetric
arrowhead matrix $A$, an let $\widetilde \lambda_{\max}$ be the eigenvalue
computed by bisection as implemented in Algorithm 3. The error bound from 
\cite[Section 3.1]{OlSt90} immediately implies that 
\begin{equation}  \label{lambda_max}
\frac{\left\vert \widetilde{\lambda }_{\max }-\lambda _{\max }\right\vert }{%
\left\vert \lambda _{\max }\right\vert }= \kappa_{bis}\varepsilon
_{M},\quad \kappa_{bis}\leq 1.06n\left( \sqrt{n}+1\right).
\end{equation}

Notice that the similar error bound holds for all eigenvalues which are of
the same order of magnitude as $\left\vert \lambda _{\max }\right\vert $.

\subsection{Accuracy of exterior eigenvalues of $A_{i}^{-1}$%
}

The desired interior eigenvalue and, in some cases, also absolutely smaller
exterior eigenvalue $\lambda$ of $A$ is in Algorithm 1 computed by (\ref{8}%
), where $\nu$ is one of the exterior eigenvalues of the matrix $A_{i}^{-1}$.

The following theorem covers the case when $\nu$ is the absolutely largest
eigenvalue of $\left\Vert A_{i}^{-1}\right\Vert _{2}$, and gives two
different bounds.

\begin{theorem}
\label{T4} Let $A_{i}^{-1}$ be defined by (\ref{invAi}) and let $\nu $ be
its eigenvalue such that 
\begin{equation}
|\nu |=\left\Vert A_{i}^{-1}\right\Vert _{2}.  \label{numax}
\end{equation}%
Let $\widehat{\nu }$ be the exact eigenvalue of the computed matrix $%
\widetilde{\left( A_{i}^{-1}\right) }=fl\left(A_{i}^{-1}\right)$. Let 
\begin{equation}
\widehat{\nu }=\nu \left( 1+\kappa _{\nu }\varepsilon _{M}\right) .
\label{kapa_nu}
\end{equation}%
Then%
\begin{equation}
\left\vert \kappa _{\nu }\right\vert \leq 
\min\bigg\{(n+3)\sqrt{n}K_{b},
 3\sqrt{n}+(n+3)\big(1+\frac{2}
{|\zeta _{i}|}\sum\limits_{\substack{ k=1  \\ k\neq i}}^{n-1}|\zeta
_{k}| \big)\bigg\},
\label{knu2}
\end{equation}%
where $K_{b}$ is defined by (\ref{K2}).
\end{theorem}

\subsection{Final error bounds}

All previous error bounds are summarized as follows.

\begin{theorem}
\label{T7} 
Let $\widetilde{\lambda }$ be the computed eigenvalue of an
unreduced arrowhead matrix $A$, let $\widetilde{\mu }$ be computed
eigenvalue of the matrix $\widetilde{A_{i}}$ from (\ref{Ai}), and let $%
\widetilde{\nu }$ be the corresponding computed eigenvalue of the matrix $%
\widetilde{\left( A_{i}^{-1}\right) }$ from (\ref{invAi}). If $\mu $ is the
eigenvalue of $A_{i}$ closest to zero (or, equivalently, if (\ref{numax})
holds), then the error in the computed eigenvalue $\widetilde{\lambda }$ is
given by (\ref{kapa_lambda}) with 
\begin{equation}
|\kappa _{\lambda }|\leq 3(|\kappa _{\nu }|+\kappa _{bis})+4,
\label{kapa_lambda_2}
\end{equation}%
and the error in the computed un-normalized eigenvector $\widetilde{x}$ is
given by Theorem \ref{T2} with 
\begin{equation}
|\kappa _{\mu }|\leq |\kappa _{\nu }|+\kappa _{bis}+1,  \label{kapa_mi_2}
\end{equation}%
where $|\kappa _{\nu }|$ is bounded by (\ref{knu2}) 
and $\kappa _{bis}$ is defined by (\ref{lambda_max}).
\end{theorem}

Since we are essentially using the shift-and-invert technique,
we can guarantee high relative accuracy of the computed
eigenvalue and high componentwise relative accuracy of the computed
eigenvector if $\nu$ is such that $|\nu|=O(\|A_{i}^{-1}\|_2)$ and it is
computed accurately.
This is certainly fulfilled if the following conditions are met:
\textit{
\begin{enumerate}
\item[C1.] The quantity $K_\nu$ from (\ref{K0}) is moderate, and
\item[C2.] 
\begin{enumerate}
\item[(i)] either the quantity $K_b$ from (\ref{K2}) is small, or
\item[(ii)] the quantity $\frac{1}{%
|\zeta _{i}|}\sum\limits_{\substack{ k=1  \\ k\neq i}}^{n-1}|\zeta
_{k}|$ from (\ref{knu2}) is of order $O(n)$.
\end{enumerate}
\end{enumerate}
}
The condition \textit{C1} implies that $\nu$ will be computed accurately
according to the standard perturbation theory.
The conditions \textit{C2 (i)} or \textit{C2 (ii)} imply that $\kappa_\nu$ from 
(\ref{knu2}) 
is small, which, together with \textit{C1}, 
implies that $\nu$ is computed accurately.

If the condition \textit{C1} does not hold, that is, if $K_\nu\gg 1$, remedies
are given in Remark 2 below.  
If neither of the conditions \textit{C2 (i)} and \textit{C2 (ii)} holds, the
remedy is to compute $b$ in double of the working precision as described in
Section 4.

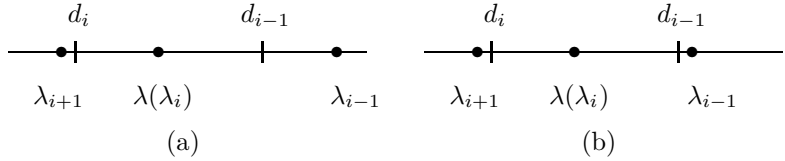
\begin{figure}[hbtp]
\begin{tabular}{cc}
\setlength{\unitlength}{0.00087489in} \begingroup\makeatletter\ifx\SetFigFont%
\undefined\gdef\SetFigFont#1#2#3#4#5{\ \reset@font\fontsize{#1}{#2pt} %
\fontfamily{#3}\fontseries{#4}\fontshape{#5} \selectfont}\fi%
\endgroup{\renewcommand{\dashlinestretch}{30}
\begin{picture}(2154,852)(0,-10)
\put(902,390){\circle*{60}}
\put(1977,390){\circle*{60}}
\put(320,390){\circle*{60}}
\drawline(2,390)(2152,390)
\drawline(402,455)(402,325)
\drawline(1527,455)(1527,325)
\put(1400,570){\makebox(0,0)[lb]{\smash{{\SetFigFont{10}{16.8}{\rmdefault}
{\mddefault}{\updefault}$d_{i-1}$}}}}
\put(1935,80){\makebox(0,0)[lb]{\smash{{\SetFigFont{10}{16.8}{\rmdefault}
{\mddefault}{\updefault}$\lambda_{i-1}$}}}} \put(150,80){\makebox(0,0)[lb]{\smash{{\SetFigFont{10}{16.8}{\rmdefault}
{\mddefault}{\updefault}$\lambda_{i+1}$}}}} \put(357,570){\makebox(0,0)[lb]{\smash{{\SetFigFont{10}{16.8}{\rmdefault}
{\mddefault}{\updefault}$d_{i}$}}}} \put(750,80){\makebox(0,0)[lb]{\smash{{\SetFigFont{10}{16.8}{\rmdefault}
{\mddefault}{\updefault}$\lambda (\lambda_i)$}}}} \end{picture} } & %
\setlength{\unitlength}{0.00087489in} \begingroup\makeatletter\ifx\SetFigFont%
\undefined\gdef\SetFigFont#1#2#3#4#5{\ \reset@font\fontsize{#1}{#2pt} %
\fontfamily{#3}\fontseries{#4}\fontshape{#5} \selectfont}\fi%
\endgroup{\renewcommand{\dashlinestretch}{30}
\begin{picture}(2154,852)(0,-10)
\put(902,390){\circle*{60}}
\put(1610,390){\circle*{60}}
\put(320,390){\circle*{60}}
\drawline(2,390)(2152,390)
\drawline(402,455)(402,325)
\drawline(1527,455)(1527,325)
\put(1400,570){\makebox(0,0)[lb]{\smash{{\SetFigFont{10}{16.8}{\rmdefault}
{\mddefault}{\updefault}$d_{i-1}$}}}}
\put(1580,80){\makebox(0,0)[lb]{\smash{{\SetFigFont{10}{16.8}{\rmdefault}
{\mddefault}{\updefault}$\lambda_{i-1}$}}}} \put(150,80){\makebox(0,0)[lb]{\smash{{\SetFigFont{10}{16.8}{\rmdefault}
{\mddefault}{\updefault}$\lambda_{i+1}$}}}} \put(357,570){\makebox(0,0)[lb]{\smash{{\SetFigFont{10}{16.8}{\rmdefault}
{\mddefault}{\updefault}$d_{i}$}}}} \put(750,80){\makebox(0,0)[lb]{\smash{{\SetFigFont{10}{16.8}{\rmdefault}
{\mddefault}{\updefault}$\lambda (\lambda_i)$}}}} \end{picture} } \\ 
(a) & (b)%
\end{tabular}
\caption{Typical situations from Remark \protect\ref{rem:2}}
\label{remedy1}
\end{figure}

\begin{remark}
\label{rem:2} 
We have two possibilities:

\begin{enumerate}
\item[(a)] we can compute $\lambda$ by shifting to another neighboring pole
provided that $K_\nu$ is in this case small (shifting to the pole $d_{i-1}$
instead of $d_i$ in Figure \ref{remedy1} (a)),

\item[(b)] if shifting to another neighboring pole is not possible ($K_\nu\gg
1$, see Figure \ref{remedy1} (b)), we can invert $A-\sigma I$, where shift $%
\sigma $ is chosen near $\lambda$, and $\sigma\notin \{\lambda,d_i,d_{i-1}\}$%
. This results in a DPR1 matrix 
\begin{equation*}
(A-\sigma I)^{-1}=%
\begin{bmatrix}
(D-\sigma I)^{-1} &  \\ 
& 0%
\end{bmatrix}%
+\rho uu^{T},  
\end{equation*}%
where 
\begin{equation*}
u=%
\begin{bmatrix}
z^{T}(D-\sigma I)^{-1} & -1%
\end{bmatrix}%
^{T},\quad \rho =\frac{1}{a-z^{T}(D-\sigma I)^{-1}z}.
\end{equation*}%
Eigenvalues of this matrix are zeros of 
\begin{equation*}
\varphi \left( \lambda \right) =1+\rho \sum_{j=1}^{n}\frac{u_{j}^{2}}{%
(d_{j}-\sigma )-\lambda },
\end{equation*}%
and the absolutely largest eigenvalue is computed accurately. Eventually, $%
\rho $ needs to be computed in higher precision.\footnote{%
Determining whether $\rho$ needs to be computed in higher precision is done
similarly as determining whether element $b$ of $A_i^{-1}$ needs to be
computed in higher precision, which is described in Section 4. Further,
Theorem \ref{T6} implies that it suffices to compute $\rho$ in double of the
working precision.}
\end{enumerate}
\end{remark}

\section{Final algorithm}

If neither of the conditions \textit{C2 (i)} and \textit{C2 (ii)} hold, in
order to guarantee that $\lambda$ will be computed with high relative
accuracy, the element $b$ from the matrix $A_{i}^{-1}$ needs to be computed
in higher precision. The following theorem implies that if $1\ll K_b\leq
O(1/\varepsilon_M)$, it is sufficient to evaluate (\ref{b_ob}) in double of the 
working precision.\footnote{If $K_b\geq O(1/\varepsilon_M)$, that is, if
  $K_b=1/\varepsilon_E$ for some $\varepsilon_E<\varepsilon_M$, then, in view of
  Theorem \ref{T6}, $b$ needs to be computed with extended precision 
$\varepsilon_E$.}
\footnote{Usage of higher precision in conjunction with the
eigenvalue computation for DPR1 matrices is analyzed in \cite{Bar93},
but there the higher precision
computation is potentially needed in the iterative part. This is less
convenient than our approach where the higher precision computation is used
only to compute one element.
}

\begin{theorem}
\label{T6} If $-a>0$ in (\ref{b_ob}), set 
\begin{equation*}
P=-a+z_{1}^{T}D_{1}^{-1}z_{1},\quad Q=-z_{2}^{T}D_{2}^{-1}z_{2},
\end{equation*}%
and if $-a<0$ in (\ref{b_ob}) set 
\begin{equation*}
P=z_{1}^{T}D_{1}^{-1}z_{1},\quad Q=a-z_{2}^{T}D_{2}^{-1}z_{2}.
\end{equation*}%
Notice that in both cases $P,Q\geq 0$ and $b=(P-Q)/\zeta _{i}^{2}$. 
Let $\widetilde P=fl(P)$ and $\widetilde Q=fl(Q)$ be evaluated in standard
precision, $\varepsilon_M$. 
Assume that $\widetilde P\neq \widetilde Q$ and $K_b\leq O(1/\varepsilon_M)$.
If 
$P$, $Q$ and $b$ are all evaluated in double of the working precision,
$\varepsilon_M^2$, then (\ref{kapa_b}) holds with $|\kappa_b|\leq O(n)$. 
\end{theorem}

We summarize the above results in one, complete algorithm, \textit{aheig}.
The algorithm first checks the components of the vector $z$. If they are of
the same order of magnitude, the eigenpair $(\lambda,v)$ is computed by
Algorithm 1. If that is not the case, the quantity $K_{b}$ is computed, and
if $K_{b}\gg 1$, the eigenpair $(\lambda,v)$ is computed by Algorithm 1 but
with evaluation of $b$ in double of the working precision. At the end, the
quantity $K_{\nu}$ is computed, and if $K_{\nu}\gg 1$, one of the remedies
from Remark \ref{rem:2} is applied.

\begin{algorithm}
\caption{}
$[\lambda,v]=\mathbf{aheig}\left(
D,z,\alpha ,k\right)$

{\small \% Computes the} $k${\small -th eigenpair
of an ordered irreducible arrowhead matrix }

{\small \%} $A=[\diag\left( D\right) \ z;
z^{\prime } \ \alpha ]$


compute the shift $i$ as in the first part of Algorithm 1

\textbf{if} {\small the quantity} $
\bigg(\sum\limits_{\substack{ j=1  \\ j\neq i}}^{n-1}|\zeta
_{j}|\bigg)/|\zeta _{i}|$ {\small from (\ref{knu2}) is of} $O(n)$ 

\ \ {\small \% standard precision is enough}

\ \ $[\lambda,v]=\mathbf{aheig\_basic}(D,z,\alpha ,k)$

\textbf{else}

\ \ compute the quantity $K_{b}$ from $(\ref{K2})$

\ \ \textbf{if} $K_{b}\gg 1$

\ \ \ \ {\small \% double precision is necessary}

\ \ \ \ $[\lambda,v]=\mathbf{aheig\_basic}(D,z,\alpha ,k)$ with evaluation of
$b$ in double precision

\ \ \textbf{else}

\ \ \ \ {\small \% standard precision is enough}

\ \ \ \ $[\lambda,v]=\mathbf{aheig\_basic}(D,z,\alpha ,k)$

\ \ \textbf{end}

\textbf{end}

compute the quantity $K_{\nu}$ from $(\ref{K0})$

\textbf{if} $K_{\nu} \gg 1$

\ \ apply one of the remedies from Remark \ref{rem:2}

\textbf{end}

\end{algorithm}

\subsection{On implementing double precision}

Implementation of the double of the working precision depends upon whether
the input is considered to be binary or decimal.

Double standard precision in Matlab, which assumes that input is binary, is
obtained by using a combination of commands \verb|vpa|, \verb|digits| and 
\verb|double| \cite{Mat}, where

\begin{itemize}
\item[-] \verb|digits(d)| specifies the number of significant decimal digits 
$d$\ used to do variable precision arithmetic \verb|vpa|,

\item[-] \verb|vpa(x)| uses variable-precision arithmetic to compute $x$ to $%
d$ decimal digits of accuracy,

\item[-] \verb|double(x)| converts $x$ to standard precision.
\end{itemize}

The assignment \verb|a1=vpa(a,32)| pads the binary representation of $a$
with zeros, which means that the decimal interpretation of the variable 
\verb|a1| may have non-zero entries after $16$-th significant decimal digit.
The same effect is obtained in Intel FORTRAN compiler \verb|ifort| \cite{Int}
by the following program segment
\begin{verbatim}
real(8) a
real(16) a1
...
a1=a
\end{verbatim}

However, the user can assume that the true input is given as a decimal
number, which is, for example, assumed by extended precision computation in
Mathematica \cite{Wol}. In this case, the options in Matlab are to either
use symbolic computation, or to cast the input to a string, and then convert
it to extended precision:
\begin{verbatim}
a1=vpa(num2str(a,16),32)
\end{verbatim}

In this case, the the decimal interpretation of the variable \verb|a1| has
all zero entries after $16$-th significant decimal digit, but the binary
representation of the variable \verb|a| is, in general, padded with non-zero
entries. The same effect is obtained in \verb|ifort| writing to and reading
from a string variable as in the following program segment:
\begin{verbatim}
real(8) a
real(16) a1
character(25) string
...
write(string,*) a
read(string,*) a1
\end{verbatim}

If the input consists of numbers for which decimal and binary
representation are equal (for example, integers, as in Example \ref{ex:3} 
below), then the two above approaches give the same results.


%

\section{Numerical Examples}

We illustrate out algorithm with four numerically demanding examples.
Examples \ref{ex:1} and \ref{ex:2} illustrate Algorithm 1, Example \ref{ex:3} 
illustrates the use of double precision arithmetic, and Example \ref{ex:quant} 
illustrates an application of higher dimension.



\begin{example}\label{ex:1}
In this example both quantities $K_\nu$ from (\ref{K0}) and $K_b$ from (\ref%
{K2}) are for all eigenvalues approximately equal to $1$, so we guarantee
that all eigenvalues and all components of their corresponding eigenvectors
are computed with high relative accuracy by Algorithm 5
(\textit{aheig}) using only standard machine precision. Let  
\begin{equation*}
A=\left[ 
\begin{array}{cccccc}
2\cdot 10^{-3} & 0 & 0 & 0 & 0 & 10^{7} \\ 
0 & 10^{-7} & 0 & 0 & 0 & 10^{7} \\ 
0 & 0 & 0 & 0 & 0 & 1 \\ 
0 & 0 & 0 & -10^{-7} & 0 & 10^{7} \\ 
0 & 0 & 0 & 0 & -2\cdot 10^{-3} & 10^{7} \\ 
10^{7} & 10^{7} & 1 & 10^{7} & 10^{7} & 10^{20}%
\end{array}%
\right] .
\end{equation*}
The eigenvalues computed by Matlab \cite{Mat} routine \verb|eig|, Algorithm
5 and Mathematica \cite{Wol} with 100 digits
precision, are, respectively: {\small 
\begin{equation*}
\begin{array}{ccc}
\lambda^{(eig)} & \lambda^{(aheig)} & \lambda^{(Math)} \\ 
\ \ 1.000000000000000\cdot 10^{20} & \ \ 1.000000000000000\cdot 10^{20} & \
\ 1.000000000000000\cdot 10^{20} \\ 
\ \ 1.999001249000113\cdot 10^{-3} & \ \ 1.999001249000113\cdot 10^{-3} & \
\ 1.999001249000113\cdot 10^{-3} \\ 
\ \ 4.987562099695390\cdot 10^{-9} & \ \ 4.987562099722817\cdot 10^{-9} & \
\ 4.987562099722817\cdot 10^{-9} \\ 
-1.000644853973479\cdot 10^{-20} & -9.999999999980001\cdot 10^{-20} & 
-9.999999999980001\cdot 10^{-20} \\ 
-2.004985562101759\cdot 10^{-6} & -2.004985562101717\cdot 10^{-6} & 
-2.004985562101717\cdot 10^{-6} \\ 
-2.001001251000111\cdot 10^{-3} & -2.001001251000111\cdot 10^{-3} & 
-2.001001251000111\cdot 10^{-3}%
\end{array}%
\end{equation*}%
}We see that even the tiniest eigenvalues $\lambda_{3}$ and $\lambda_{4}$,
computed by Algorithm 5, are exact to the machine precision, which
is not true for the eigenvalues computed by \verb|eig|. Because of the
accuracy of the computed eigenvalues, the eigenvectors computed by 
Algorithm 5 are componentwise accurate up to machine precision, and
therefore, orthogonal up to machine precision. For example: 
{\small 
\begin{equation*}
\begin{array}{ccc}
v_{4}^{(eig)} & v_{4}^{(aheig)} & v_{4}^{(Math)} \\ 
\ \ 4.999993626151683\cdot 10^{-11} & -4.999999999985000\cdot 10^{-11} & 
-4.999999999985000\cdot 10^{-11} \\ 
\ \ 9.999999962328609\cdot 10^{-7} & -9.999999999969000\cdot 10^{-7} & 
-9.999999999969000\cdot 10^{-7} \\ 
\ \ 9.999999999990000\cdot 10^{-1} & -9.999999999989999\cdot 10^{-1} & 
-9.999999999989999\cdot 10^{-1} \\ 
-9.999999964673912\cdot 10^{-7} & \ \ 9.999999999970999\cdot 10^{-7} & \text{
\ }9.999999999970999\cdot 10^{-7} \\ 
-5.000006338012225\cdot 10^{-11} & \ \ 4.999999999985000\cdot 10^{-11} & \ \
4.999999999985000\cdot 10^{-11} \\ 
-9.999999963825105\cdot 10^{-21} & \ \ 9.999999999970000\cdot 10^{-21} & \ \
\ 9.999999999969999\cdot 10^{-21}.%
\end{array}%
\end{equation*}%
} 
\end{example}


\begin{example}\label{ex:2}
In this example, despite very close diagonal elements, we again guarantee
that all eigenvalues and all components of their corresponding eigenvectors
are computed with high relative accuracy, without deflation. Let 
\begin{equation*}
A=\left[ 
\begin{array}{ccccc}
1+4\varepsilon _{M} & 0 & 0 & 0 & 1 \\ 
0 & 1+3\varepsilon _{M} & 0 & 0 & 2 \\ 
0 & 0 & 1+2\varepsilon _{M} & 0 & 3 \\ 
0 & 0 & 0 & 1+1\varepsilon _{M} & 4 \\ 
1 & 2 & 3 & 4 & 0%
\end{array}%
\right]
\end{equation*}%
where $\varepsilon _{M}=2\cdot 2^{-53}=2.2204\cdot 10^{-16}$. 
For this matrix the quantities $K_{\nu} $ and $K_{b} $ are again of order one
for all eigenvalues,
so Algorithm 5 uses only standard working precision. 
The eigenvalues computed by Matlab and Algorithm 5 are:
{\small 
\begin{equation*}
\begin{array}{cc}
\lambda ^{(eig)} & \lambda ^{(aheig)} \\ 
6.000000000000000 & 6.000000000000001 \\ 
1+4\varepsilon _{M} & 1+4\varepsilon _{M} \\ 
1+4\varepsilon _{M} & 1+3\varepsilon _{M} \\ 
1+3\varepsilon _{M} & 1+2\varepsilon _{M} \\ 
-5.000000000000000 & -4.999999999999999 
\end{array}%
\end{equation*}
}
The eigenvalues computed by 
Mathematica with 100 digits precision, properly rounded to 32 decimal digits 
are\footnote{Since, as described in Section 4.1,
  Mathematica uses decimal representation of the input, in order to obtain
  accurate eigenvectors we need to define $\varepsilon_M$ in Mathematica with 
  the output of Matlab's command vpa(eps), 
  $\varepsilon_M=2.2204460492503130808472633361816\cdot 10^{-16}$.}:
{\small 
\begin{equation*}
\begin{array}{c}
\lambda ^{(Math)} \\ 
6.0000000000000002018587317500285 \\ 
1.0000000000000008727792604471857 \\ 
1.0000000000000006206061701073114\\ 
 1.0000000000000003571862771540971\\ 
-4.9999999999999998317843902083010 
\end{array}%
\end{equation*}
}The eigenvalues computed by Matlab are accurate according to standard
perturbation theory, but they do not satisfy the interlacing property. 
Furthermore, the Matlab's eigenvectors corresponding to $\lambda_2$,
$\lambda_3$ and $\lambda_4$ only span an accurate eigenspace, and are not
individually accurate. On the other hand, 
the eigenvalues computed by Algorithm 5 are exact (they coincide
with the eigenvalues computed by Mathematica properly rounded to 16 decimal
digits). Notice that despite of very close eigenvalues, Algorithm 5 works
without deflation. Due to the accuracy of the computed
eigenvalues, the eigenvectors computed by Algorithm 5 are
componentwise accurate up to the machine precision, and are therefore 
orthogonal.

If, as suggested in Remark \ref{rem:1a}, the algorithms are modified to return
$d_i$ and $\mu$ (both in standard precision), then for the eigenvalues
$\lambda_2$, $\lambda_3$ and $\lambda_4$ the corresponding pairs $(d_i,\mu)$
give representations of those eigenvalues to 32 decimal digits. That is, 
exact values $d_i+\mu$ properly rounded to 32 decimal digits are equal to 
the corresponding eigenvalues computed by Mathematica as displayed above.


\end{example}


\begin{example}\label{ex:3}
In this example we can guarantee all eigenvalues and eigenvectors,
componentwise will be computed with high relative accuracy only if we use
double of the working precision when computing $b$ from (\ref{b_ob}) 
in matrices $%
A_{2}^{-1},A_{3}^{-1},A_{4}^{-1}$ and $A_{5}^{-1}$. Let 
\begin{equation}
A=\left[ 
\begin{array}{cccccc}
10^{10} & 0 & 0 & 0 & 0 & 10^{10} \\ 
0 & 4 & 0 & 0 & 0 & 1 \\ 
0 & 0 & 3 & 0 & 0 & 1 \\ 
0 & 0 & 0 & 2 & 0 & 1 \\ 
0 & 0 & 0 & 0 & 1 & 1 \\ 
10^{10} & 1 & 1 & 1 & 1 & 10^{10}%
\end{array}%
\right] .  \notag
\end{equation}

The quantities $K_{\nu} $ and $K_{b} $ are:\footnote{Algorithm 5 does not
  compute $K_{\nu} $ and $K_{b}$ for the first eigenvalue, since it is an
  absolutely largest one.}
{\small 
\begin{equation*}
\begin{array}{cc}
K_{\nu}  & K_{b} \\ 
- & - \\ 
9.999999090793056\cdot 10^{-1} & 3.243243243540540\cdot 10^{9} \\ 
9.999996083428923\cdot 10^{-1} & 3.636363636818182\cdot 10^{9} \\ 
1.000000117045544\cdot 10^{0} & 4.444444445000000\cdot 10^{9} \\ 
9.999998561319470\cdot 10^{-1} & 5.217390439488477\cdot 10^{9} \\ 
7.941165469988994\cdot 10^{0} & 5.217390439488477\cdot 10^{9}%
\end{array}%
\end{equation*}%
}It is clear, from the condition numbers, that the element $b$ in each of the
matrices $%
A_{2}^{-1},A_{3}^{-1},A_{4}^{-1}$ and $A_{5}^{-1}$ needs to be computed in
double of the working precision.
For example,%
\begin{equation*}
A_{2}=A-d_{2}I=\left[ 
\begin{array}{cccccc}
10^{10}-4 & 0 & 0 & 0 & 0 & 10^{10} \\ 
0 & 0 & 0 & 0 & 0 & 1 \\ 
0 & 0 & -1 & 0 & 0 & 1 \\ 
0 & 0 & 0 & -2 & 0 & 1 \\ 
0 & 0 & 0 & 0 & -3 & 1 \\ 
10^{10} & 1 & 1 & 1 & 1 & 10^{10}-4%
\end{array}%
\right] .
\end{equation*}%
The element $b=\left[ A_2^{-1}\right]_{22}$ computed by Algorithm 2
gives $b=6.16666666667$, Matlab's routine \verb!inv! yields
$b=6.166665889418350$, while computing $b$ in double of the working precision
gives the correct value $b=6.166666668266667$.


Eigenvalues computed by Algorithm 1 (\textit{aheig\_basic}, using only
standard working precision), Algorithm 5 (\textit{aheig}, using double of the
working precision to compute respective $b$'s) and
Mathematica with 100 digits precision, respectively, are: {\small 
\begin{equation*}
\begin{array}{ccc}
\lambda ^{_{aheig\_basic}} & \lambda ^{_{^{aheig}}} & \lambda ^{_{Math}} \\ 
\ \ 2.000000000000000\cdot 10^{10} & \text{ \ }2.000000000000000\cdot 10^{10}
& \ \ 2.000000000000000\cdot 10^{10} \\ 
\ \ 4.150396802313551\cdot 10^{0} & \ \ 4.150396802279712\cdot 10^{0} & \ \
4.150396802279713\cdot 10^{0} \\ 
\ \ 3.161498641452035\cdot 10^{0} & \ \ 3.161498641430967\cdot 10^{0} & \ \
3.161498641430967\cdot 10^{0} \\ 
\ \ 2.188045596352105\cdot 10^{0} & \ \ 2.188045596339914\cdot 10^{0} & \ \
2.188045596339914\cdot 10^{0} \\ 
\ \ 1.216093560005649\cdot 10^{0} & \ \ 1.216093584948579\cdot 10^{0} & \ \
1.216093584948579\cdot 10^{0} \\ 
-7.160348702977373\cdot 10^{-1} & -7.160346250991725\cdot 10^{-1} & 
-7.160346250991725\cdot 10^{-1}%
\end{array}%
\end{equation*}%
}The eigenvectors computed by Algorithm 5 are componentwise
accurate to machine precision and therefore orthogonal.
\end{example}


\begin{example}
\label{ex:quant}
This example comes from the research 
related to decay of excited states of quantum dots in in real photon crystals
\cite{MM10}. In this case 
\begin{itemize}
\item[-] $\alpha$ is quantum dot transition frequency,

\item[-] $d_{i}$ is a frequency of the $i$-th optical mode, and

\item[-] $\zeta_{i}$ is an interaction constant of the quantum dot with the
  $i$-th optical mode.
\end{itemize}
The size of the matrix is changeable but, in realistic cases, it is
between $10^{3}$ and $10^{4}$. We ran a test example for 
$n=2501$ where, typically, 
\begin{align*}
d_i & \in [5.87 \cdot 10^{14},1.38 \cdot 10^{15}],\\
  \zeta_i& \in [1.05\cdot 10^{4},1.10 \cdot 10^{7}],\\
\alpha&=9.7949881500060375\cdot 10^{14}.
\end{align*}%
For this matrix the condition number $K_{\nu}\sim 1$ for all
eigenvalues and  the components of the vector $z$ do not differ by much in
size, thus the conditions \textit{C1} and \textit{C2 (ii)} from Section
3 are fulfilled. Therefore, 
all eigenvalues and all components of all eigenvectors
are computed with high relative accuracy by Algorithm 
5 using only standard working precision. On the other hand about half of
the eigenvalues computed by the Matlab routine \verb!eig! do not satisfy the
interlacing property.

\end{example}




\section{Applications}
In this  section we extend our results to 
eigenvalue decompositions of Hermitian arrowhead matrices, singular value
decompositions of real triangular arrowhead matrices and  
eigenvalue decompositions of real symmetric diagonal-plus-rank-one matrices. 

\subsection{Hermitian arrowhead matrices}


Let

\begin{equation*}
C=\left[ 
\begin{array}{cc}
D & z \\ 
z^{*} & \alpha%
\end{array}%
\right],
\end{equation*}%
where 
\begin{equation*}
D=\mathop{\mathrm{diag}}(d_{1},d_{2},\ldots ,d_{n-1}),
\end{equation*}%
is a real diagonal matrix of order $n-1$,%
\begin{equation*}
z=\left[ 
\begin{array}{cccc}
\zeta _{1} & \zeta _{2} & \cdots & \zeta _{n-1}%
\end{array}%
\right] ^{*},
\end{equation*}%
is a complex valued vector and $\alpha$ is a real scalar. Here $z^*$ denotes
the conjugate transpose of $z$. As in Section 1, we assume that $C$ is
irreducible. The eigenvalue decomposition of $C$ is given by 
\begin{equation*}
C=U\Lambda U^{\ast }
\end{equation*}%
where $\Lambda=\mathop{\mathrm{diag}}(\lambda_1,\ldots,\lambda_n)\in\mathbb{R%
}^{n\times n}$ is a diagonal matrix of eigenvalues, and $U=%
\begin{bmatrix}
u_1 & u_2 & \cdots & u_n%
\end{bmatrix}%
$ is an unitary matrix of the corresponding eigenvectors.

To apply Algorithm 5 to Hermitian arrowhead matrix we first transform $C$ to
real symmetric arrowhead matrix $A$ by diagonal unitary similarity:

\begin{equation}  \label{cc}
A=\Phi^*C \Phi=\left[ 
\begin{array}{cc}
D & |z| \\ 
|z|^{T} & \alpha%
\end{array}%
\right] ,
\end{equation}
where 
\begin{equation*}
\Phi=\mathop{\mathrm{diag}}\left( \frac{\zeta_{1}}{\left\vert
\zeta_{1}\right\vert }, \frac{\zeta_{2}}{\left\vert \zeta_{2}\right\vert }%
,\ldots, \frac{\zeta_{n-1}}{\left\vert \zeta_{n-1}\right\vert },1\right)
\end{equation*}

We now compute the $k$-th eigenpair $(\lambda,v)$ of $A$ by Algorithm 5, and
set $u=\Phi v$. Since we guarantee high relative accuracy of the eigenvalue
decomposition of $A$ computed by Algorithm 5, we also guarantee high
relative accuracy of the eigenvalue decomposition of $C$. Notice that, if
double precision is needed to compute $b$ in Algorithm 5, the modules $%
|\zeta_i|$ in (\ref{cc}) need to be computed in double of the working
precision, as well.

\begin{remark}
Similarly, for irreducible non-symmetric arrowhead matrix 
\begin{equation*}
G=\left[ 
\begin{array}{cc}
D & \breve{z} \\ 
\mathring{z}^T & \alpha%
\end{array}%
\right],
\end{equation*}%
where $\sign(\mathring\zeta_i)=\sign (\breve\zeta_i)$, $i=1,\ldots,n-1$,
we define the diagonal matrix 
\begin{equation*}
\Psi=\diag\left( \sign
  (\mathring\zeta_{1})\sqrt{\frac{\mathring\zeta_{1}}{\breve\zeta_{1}}},\ldots, 
\sign
  (\mathring\zeta_{n-1})\sqrt{\frac{\mathring\zeta_{n-1}}{\breve\zeta_{n-1}}},1\right).
\end{equation*}
The matrix
\begin{equation*}
A=\Psi^{-1}G \Psi=\left[ 
\begin{array}{cc}
D & z \\ 
z^{T} & \alpha%
\end{array}%
\right] ,
\end{equation*}
where $\zeta_i=\sqrt{\breve\zeta_{1} \mathring\zeta_i}$ is an irreducible
symmetric arrowhead matrix.

We now compute the $k$-th eigenpair $(\lambda,v)$ of $A$ by Algorithm 5. The
eigenpair of $G$ is then $(\lambda,\Psi v)$.
set $u=\Phi v$. Since we guarantee high relative accuracy of the eigenvalue
decomposition of $A$, we also guarantee high
relative accuracy of the eigenvalue decomposition of $G$. Notice that, if
double precision is needed to compute $b$ in Algorithm 5, the elements
$\zeta_i$ need to be computed in double of the working
precision, as well.
\end{remark}

\subsection{Singular value decomposition of a triangular arrowhead matrix}

Let 
\begin{equation*}
B=\left[ 
\begin{array}{cc}
D & z \\ 
0  & \alpha%
\end{array}%
\right],
\end{equation*}%
be an irreducible upper triangular arrowhead matrix, that is, $d_i\neq d_j$ for
$i\neq j$ and $\zeta_i\neq 0$ for all $i$.
The matrix 
\begin{equation*}
A=B^TB=\left[ 
\begin{array}{cc}
D^2 & Dz \\ 
z^T D   & \alpha+z^T z%
\end{array}%
\right],
\end{equation*}%
is an irreducible symmetric arrowhead matrix.

When applying Algorithm 5 to the matrix $A$, we must ensure that all
components of $A_i^{-1}$
in (\ref{invAi}) are computed to high relative accuracy. This is obviously
true for elements of the vectors $w_i$ and $w_2$. Diagonal elements, except
$b$, are computed with high relative accuracy as differences of squares of 
original quantities, 
$$
[A_i^{-1}]_{jj}=\frac{1}{(d_j-d_i)(d_j+d_i)}, \quad j\neq i.
$$
The element $b=[A_i^{-1}]_{ii}$ from (\ref{b_ob}) is computed as\footnote{In
  view of Theorem \ref{T6}, if double precision computation is necessary, the
  positive and negative parts of this formula should be computed separately,
  and then added.}
$$
b=\frac{1}{d_i^2 \zeta_i^2}\left(-\alpha-z^Tz +d_i+\sum_{j \neq i}
\frac{d_j^2\zeta_j^2}{(d_j-d_i)(d_j+d_i)}\right).
$$
If double precision is needed in Algorithm 5, all 
entries of $A$ need to be computed in double precision.

Let $B=U\Sigma V^T$ be the singular value decomposition of $B$, where 
$\Sigma=\diag(\sigma_1,\ldots,\sigma_n)$ are the singular values, 
the columns of $V$ are the corresponding right singular vectors and 
the columns of $U$ are the corresponding left singular vectors.
%
We first compute the $k$-th eigenpair $(\lambda,v)$ of $A$ by Algorithm 5.
Then $\sigma=\sqrt{\lambda}$ is the corresponding singular value of $B$ and 
$v$ is the corresponding right singular vector. The value $\sigma$ and all
components of $v$ are computed to almost full accuracy.
From the relation $U^T B=\Sigma V^T$ for the $k$-th row we have
$$
\begin{bmatrix} u_{1:n-1}^T & u_n
\end{bmatrix}
\begin{bmatrix} D & z \\ 0 & \alpha
\end{bmatrix}=\sigma \begin{bmatrix} v_{1:n-1}^T & v_n
\end{bmatrix},
$$
which implies
$$
u_{1:n-1}= \sigma v_{1:n-1} D^{-1}.
$$
From the relation $BV=U\Sigma $ for the $k$-th column we have
$$
\begin{bmatrix} D & z \\ 0 & \alpha
\end{bmatrix}
\begin{bmatrix} v_{1:n-1} \\ v_n
\end{bmatrix}
=\sigma \begin{bmatrix} u_{1:n-1} \\ u_n
\end{bmatrix},
$$
which implies
$$
u_n= \frac{\alpha v_n}{\sigma}.
$$
Components of $u$ are computed by multiplication and division of
quantities which are accurate to almost full machine precision, so the are
accurate to almost full machine precision, as well.

\subsection{Diagonal-plus-rank-one matrices}

Let
$$
M =D+u u^T,
$$
where
\begin{align*}
D&=\diag(d_1,\ldots,d_n), \quad d_1>d_2>\cdots>d_n,\\ 
u&=\begin{bmatrix}u_1 & \cdots & u_n \end{bmatrix}^T, 
\quad u_i\neq 0, \quad i=1,\ldots,n,
\end{align*}
be a $n\times n$ irreducible ordered real symmetric diagonal-plus-rank-one
(DPR1) matrix. Let 
\begin{align*}
\bar D& =\diag(d_1,\ldots,d_{n-1}), \\
\Delta&=(\bar D-d_n)^{1/2},\\
\bar u&=\begin{bmatrix}u_1 & \cdots & u_{n-1} \end{bmatrix}^T,\\
L&=\begin{bmatrix}
u_n \Delta^{-1}& 0 \\
-\bar u^T \Delta^{-1} & 1 
\end{bmatrix}.
\end{align*}
Then 
$$
A=L^{-1} M L =
\begin{bmatrix} \bar D & z\\z^T & \alpha
\end{bmatrix},
$$
where
$$
\quad z=\bar u \Delta, \quad \alpha= d_{n}+u^T u,
$$
is an irreducible real symmetric arrowhead matrix.

When applying Algorithm 5 to the matrix $A$, we must ensure that all
components of $A_i^{-1}$
in (\ref{invAi}) are computed to high relative accuracy. This is obviously
true for elements of the vectors $w_i$ and $w_2$. Diagonal elements, except
$b$, are computed with high relative accuracy as differences of 
original quantities, and the element
$b=[A_i^{-1}]_{ii}$ from (\ref{b_ob}) is
computed as
$$
b=\frac{1}{\zeta_i}\left(-d_n-u^Tu +d_i+\sum_{j \neq i}
\frac{\zeta_j}{d_j-d_i}\right).
$$
If double precision is needed in Algorithm 5, all 
entries of $A$ need to be computed in double precision.

Let $M=Q\Lambda Q^T$ and $A=V\Lambda V^T$ be the eigenvalue decompositions of 
$M$ and $A$, respectively. Since $M$ is by assumption irreducible, its
eigenvalues satisfy interlacing property
\begin{equation}\label{lld}
\lambda_1> d_1>\lambda_2>d_2> \cdots> \lambda_n>d_n.
\end{equation}
We first compute the $k$-th eigenpair $(\lambda,v)$ of $A$ by Algorithm 5.
The value $\lambda$ and all
components of $v$ are computed to almost full accuracy.
The relation $V^T AV=V^T L^{-1} M L V= \Lambda $ implies that the columns of
the matrix $X=LV$ are the unnormalized eigenvectors of the matrix $M$. Further,
since, by (\ref{lld}), all eigenvalues are simple, we conclude that 
$X=Q\Sigma$, where $\Sigma=\diag(\sigma_1,\ldots,\sigma_n)$ is a positive
definite matrix. Notice that $Q\Sigma V^T=L$ is, in fact, singular value
decomposition of $L$.

Equating $k$-th columns of the equation $X=LV$ gives
$$
x=\begin{bmatrix} \bar x\\ x_n
\end{bmatrix}
=Lv=\begin{bmatrix}
u_n \Delta^{-1}& 0 \\
-\bar u^T\Delta^{-1} & 1 
\end{bmatrix}
\begin{bmatrix} \bar v\\ v_n
\end{bmatrix},
$$ 
where $x$ and $v$ are partitioned according to $L$.
This immediately implies that
$$
\bar x=u_n \Delta^{-1} \bar v.
$$
Notice that, since all components of $\bar v$ are computed to almost full,
accuracy, the same holds for the components of $\bar x$, and it remains to
compute $x_n$ accurately.

Let 
$$
q=\begin{bmatrix} \bar q \\ q_n
\end{bmatrix}
$$
be the $k$-th column of $Q$ and let $\sigma=\Sigma_{kk}$.
Equating $k$-th rows of the equation
$$
X^{-1}=\Sigma^{-1}Q^T= V^T L^{-1}
$$
gives for the $n$-th element 
$$
q_n\frac{1}{\sigma}= x_n\frac{1}{\sigma^2}=v_n.
$$
Thus,
$$
x_n=\sigma^2 v_n
$$
and, in order to compute $x_n$, it is necessary to compute $\sigma^2$.
From $X=U\Sigma=LV$ it follows that $V^T L^T L V=\Sigma^2$, or, equivalently, 
$LV=L^{-T} V\Sigma^2$.
Equating $k$-th columns of this equation gives
$$
\Delta^{-1} \bar v\,  u_n = \left[\Delta\bar v \frac{1}{u_n} + 
\Delta^{-1}\bar u\, v_n\right]\sigma^2.
$$
This gives $n-1$ equations for $\sigma^2$, and we can choose the
numerically most accurate one. 

Therefore, $x_n$ will be computed to almost full machine precision, as are
the entries of $\bar x$, and it remains to normalize $x$ and obtain
$q=x/\sigma$.



\begin{remark}
Notice that DPR1 matrices of the form $D-uu^T$ cannot be reduced to
symmetric arrowhead matrix by the procedure described in this section. 
By using ideas from this paper, it is possible to derive highly accurate
algorithm for DPR1 matrices without prior transformation to arrowhead
form. This algorithm, 
which is a topic of our forthcoming paper, covers more general
DPR1 matrices of the form 
$$
D+\rho u u^T, \quad \rho \in \mathbb{R}.
$$

\end{remark}

\appendix
\label{App}

\section{Proofs}

\subsection*{Proof of Theorem \ref{T1}.}

Let $\widetilde{\mu}$ and $\widetilde{\lambda}$ be defined by (\ref{tablica}).
Then
\begin{equation*}
\widetilde{\lambda }\equiv fl\left( d_{i}+\widetilde{\mu }\right) =\left(
d_{i}+\widetilde{\mu }\right) \left( 1+\varepsilon _{1}\right). 
\end{equation*}%
By simplifying the equality 
\begin{equation*}
\left( d_{i}+\mu \left( 1+\kappa _{\mu }\varepsilon _{M}\right) \right)
\left( 1+\varepsilon _{1}\right) =\lambda \left( 1+\kappa _{\lambda
}\varepsilon _{M}\right)  
\end{equation*}%
and using $\lambda =\mu +d_{i}$, we have%
\begin{equation*}
d_{i}\varepsilon _{1}+\mu \left( \kappa _{\mu }\varepsilon _{M}+\varepsilon
_{1}\right) =\lambda \kappa _{\lambda }\varepsilon _{M}.
\end{equation*}%
Taking absolute value gives%
\begin{equation*}
\left\vert \kappa _{\lambda }\right\vert \leq \frac{\left\vert
d_{i}\right\vert +\left\vert \mu \right\vert }{\left\vert \lambda
\right\vert }\left( \left\vert \kappa _{\mu }\right\vert +1\right).
\blacksquare
\end{equation*}

\subsection*{Proof of Theorem \ref{C1}.}

\textit{(i)} The assumption $\sign\left( d_{i}\right) =\sign\left( \mu
\right) $ immediately implies
\begin{equation*}
\frac{\left\vert d_{i}\right\vert +\left\vert \mu \right\vert }{\left\vert
\lambda \right\vert }=\frac{\left\vert d_{i}+\mu \right\vert }{\left\vert
d_{i}+\mu \right\vert }=1.
\end{equation*}

\textit{(ii)} The assumptions imply that either 
$$
0<d_{i+1}<\lambda <d_{i}, \quad \mu <0,
$$
or 
$$
d_{i}<\lambda <d_{i-1}<0, \quad \mu >0.
$$
In the first case $\lambda $ is closest to the pole $d_{i}$ and 
\begin{align*}
\frac{\left\vert d_{i}\right\vert +\left\vert \mu \right\vert }{\left\vert
\lambda \right\vert }& \leq \frac{\left\vert d_{i}\right\vert +\frac{1}{2}%
\left\vert d_{i}-d_{i+1}\right\vert }{\frac{1}{2}\left\vert
d_{i}+d_{i+1}\right\vert }
\leq \frac{d_{i}+\frac{1}{2}d_{i}-\frac{1}{2}d_{i+1}}{\frac{1}{2}d_{i}+%
\frac{1}{2}d_{i+1}} \\
& \leq \frac{\frac{3}{2}d_{i}-\frac{1}{2}d_{i+1}}{\frac{1}{2}d_{i}+\frac{1}{2%
}d_{i+1}}
\leq \frac{3d_{i}}{d_{i}}=3.
\end{align*}%
Here we used the inequalities $\left\vert \mu \right\vert \leq \frac{1}{2}%
\left\vert d_{i}-d_{i+1}\right\vert $ and $\left\vert \lambda \right\vert
\geq \frac{1}{2}\left\vert d_{i}+d_{i+1}\right\vert $ for the first
inequality, $d_{i}-d_{i+1}>0$ and $d_{i}+d_{i+1}>0$ for the second inequality
and $d_{i+1}>0$ for the fourth inequality, respectively.

The proof for the second case is analogous.
$\blacksquare $

\subsection*{Proof of Theorem \ref{T2}.}

Let $x$ and $\widetilde{x}$ be defined by (\ref{eigenvecAi}) and
(\ref{x_tilda}), respectively.
The theorem obviously holds for $x_{n}=\widetilde{x}_{n}=-1$.
For $\widetilde{x}_i$ we have
\begin{equation*}
\widetilde{x}_{i}=fl\left( -\frac{\zeta _{i}}{\widetilde{\mu }}\right) =-%
\frac{\zeta _{i}}{\mu \left( 1+\kappa _{\mu }\varepsilon _{M}\right) }\left(
1+\varepsilon _{1}\right) =x_{i}\left( 1+\varepsilon _{x_{i}}\right) .
\end{equation*}%
By using (\ref{kapa_mi}) and (\ref{x_tilda}), the first order approximation 
gives
\begin{equation*}
\left\vert \varepsilon _{x_{i}}\right\vert \leq \left( |\kappa _{\mu
}|+1\right) \varepsilon _{M}.
\end{equation*}%
For $j\notin\{ i,n\}$, by solving the equality 
\begin{equation*}
\widetilde{x}_{j}=\frac{\zeta _{j}}{\left( \left( d_{j}-d_{i}\right) \left(
1+\varepsilon _{1}\right) -\mu \left( 1+\kappa _{\mu }\varepsilon
_{M}\right) \right) \left( 1+\varepsilon _{2}\right) }\left( 1+\varepsilon
_{3}\right) =\frac{\zeta _{j}}{d_{j}-\lambda }\left( 1+\varepsilon
_{x}\right)
\end{equation*}%
for $\varepsilon _{x}$, using (\ref{kapa_mi}) and $\lambda =\mu +d_{i}$, and
ignoring higher order terms, we have 
\begin{equation*}
\varepsilon _{x}=\frac{\left( d_{j}-d_{i}\right) \left( \varepsilon
_{1}+\varepsilon _{2}+\varepsilon _{3}\right) -\mu \left( \kappa _{\mu
}\varepsilon _{M}+\varepsilon _{2}+\varepsilon _{3}\right) }{d_{j}-\lambda }.
\end{equation*}%
Therefore, 
\begin{equation}
\left\vert \varepsilon _{x}\right\vert \leq \frac{\left\vert
d_{j}-d_{i}\right\vert +\left\vert \mu \right\vert }{\left\vert
d_{j}-\lambda \right\vert }\left( |\kappa _{\mu }|+3\right) \varepsilon _{M}.
\label{eq:52}
\end{equation}%

To complete the proof we need to analyze two cases.
If 
$$
\mathop{\mathrm{sign}}\left( d_{j}-d_{i}\right) =-%
\mathop{\mathrm{sign}}\mu,
$$ 
then 
\begin{equation*}
\frac{\left\vert d_{j}-d_{i}\right\vert +\left\vert \mu\right\vert }{%
\left\vert d_{j}-\lambda\right\vert }=\frac{|d_{j}-d_{i}-\mu|}{\left\vert
d_{j}-\lambda\right\vert }=\frac{|d_{j}-\lambda|}{\left\vert
d_{j}-\lambda\right\vert }=1.
\end{equation*}
If 
$$
\mathop{\mathrm{sign}}\left( d_{j}-d_{i}\right) =%
\mathop{\mathrm{sign}}\mu,
$$ 
then, since $d_{i}$ is pole closest to $\lambda $,
we have $\left\vert \mu \right\vert \leq 0.5\left\vert
d_{j}-d_{i}\right\vert $ and 
\begin{equation*}
\frac{\left\vert d_{j}-d_{i}\right\vert +\left\vert \mu \right\vert }{%
\left\vert d_{j}-\lambda \right\vert }\leq \frac{\left\vert
d_{j}-d_{i}\right\vert +\left\vert \mu \right\vert }{\left\vert
d_{j}-d_{i}\right\vert -\left\vert \mu \right\vert }\leq \frac{\frac{3}{2}%
\left\vert d_{j}-d_{i}\right\vert }{\frac{1}{2}\left\vert
d_{j}-d_{i}\right\vert }=3.
\end{equation*}%
Finally, the theorem follows by  inserting this into (\ref{eq:52}).
$\blacksquare$

\subsection*{Proof of Theorem \ref{T3}.}

For the non-zero computed elements of the matrix $A_{i}^{-1}$ from 
(\ref{invAi}) and (\ref{b_ob}), except the element $b=[A_i^{-1}]_{ii}$, we have: 
\begin{align*}
fl ( \left[A_i^{-1}\right] _{jj} ) & =\frac{1}{\left(
d_{j}-d_{i}\right) \left( 1+\varepsilon_{1}\right) }\left( 1+\varepsilon
_{2}\right), \quad j\notin\{i,n\},\\
fl(\left[A_i^{-1} \right]_{ji}) & =
fl (\left[A_i^{-1} \right]_{ij})
=\frac{-\zeta _{j}}{\left(
d_{j}-d_{i}\right) \left( 1+\varepsilon_{3}\right) \zeta _{i} 
(1+\varepsilon_{4})}\left(
1+\varepsilon _{5}\right), \ j\notin\{i,n\},
\\
fl([A_i^{-1}] _{ni}) &= 
fl([A_i^{-1}] _{in})=
\frac{1}{\zeta_i}(1+\varepsilon_6),
\end{align*}%
where $|\varepsilon_k|\leq \varepsilon_M$ for all indices $k$.
The first statement of the theorem now follows by using standard first order 
approximations.

Similar analysis of the formula (\ref{b_ob}) yields
$$
fl([A_i^{-1}] _{ii})=\widetilde{b}=b+\delta b,
$$ 
where 
\begin{equation}\label{deltab}
\left\vert \delta b\right\vert \leq \frac{1}{\zeta _{i}^{2}}\left(
\left\vert -a\right\vert +\left\vert z_{1}^{T}D_{1}^{-1}z_{1}\right\vert
+\left\vert z_{2}^{T}D_{2}^{-1}z_{2}\right\vert \right)(n+3) \varepsilon _{M}.
\end{equation}
This, in turn, implies (\ref{kapa_b}) with
\begin{equation*}
\left\vert \kappa _{b}\right\vert \leq \frac{\left\vert \delta
b\right\vert }{\left\vert b\right\vert }\, \frac{1}{\left\vert
\varepsilon _{M}\right\vert }
=(n+3)\frac{\left\vert a\right\vert +\left\vert
z_{1}^{T}D_{1}^{-1}z_{1}\right\vert +\left\vert
z_{2}^{T}D_{2}^{-1}z_{2}\right\vert }{\left\vert
-a+z_{1}^{T}D_{1}^{-1}z_{1}+z_{2}^{T}D_{2}^{-1}z_{2}\right\vert }
=(n+3)K_{b},
\end{equation*}
where $K_{b}$ is defined by (\ref{K2}).
$\blacksquare$
\subsection*{Proof of Theorem \ref{T4}.}

Let
\begin{equation*}
\widetilde{\left( A_{i}^{-1}\right) }=A_{i}^{-1}+\delta A_{i}^{-1}.
\end{equation*}
Therefore,
$$
|\widehat{\nu}-\nu|=\| \delta A_{i}^{-1}\|_2,
$$
which, together with (\ref{kapa_nu}), implies 
\begin{equation}\label{kappa_nu1}
|\nu \kappa _{\nu }\varepsilon _{M}|\leq \| \delta
A_{i}^{-1}\|_{2}.
\end{equation}
Theorem \ref{T3} implies that 
$$
\| \delta A_{i}^{-1}\|_2
\leq 
(n+3)\| |A_i^{-1} |\|_2 K_b\, \varepsilon_M.
$$
Since $\| |A_i^{-1} |\|_2\leq\sqrt{n}\| A_i^{-1} \|_2$ and
$|\nu|=\|A_{i}^{-1}\|_2$  , from (\ref{kappa_nu1}) we have 
\begin{equation}\label{konacno1}
|\kappa _{\nu }|\leq (n+3)\sqrt{n}K_b,
\end{equation}
which proves the first part of the bound (\ref{knu2}).

For the second part of the proof, notice that Theorem
\ref{T3} also implies

\begin{equation}\label{nukanu}
 \| \delta
A_{i}^{-1}\|_{2}\leq 
3 \| |\mathcal{A}|\|_2 \varepsilon_M +|\delta b|, 
\end{equation}
where $\mathcal{A}$ is equal to the matrix $A_i^{-1}$ without $b$ (that is,
with $\mathcal{A}_{ii}=0$).

By bounding (\ref{kappa_nu1}) with (\ref{nukanu}) and (\ref{deltab}), and
dividing the resulting inequality by $|\nu\varepsilon_m|$, we have
\begin{equation}
|\kappa _{\nu }|\leq 3 \sqrt{n}+(n+3)\bigg(\frac{1}{\left\vert \nu
\right\vert }\frac{|-a|+\left\vert z_{1}^{T}D_{1}^{-1}z_{1}\right\vert
+\left\vert z_{2}^{T}D_{2}^{-1}z_{2}\right\vert }{\zeta _{i}^{2}}\bigg) .
\label{eq:508}
\end{equation}%
Since
\begin{align*}
\frac{|-a|}{\zeta_i^2}&=\frac{1}{\zeta_i^2}|-a+z_{1}^{T}D_{1}^{-1}z_{1}+
z_{2}^{T}D_{2}^{-1}z_{2} 
-z_{1}^{T}D_{1}^{-1}z_{1}-z_{2}^{T}D_{2}^{-1}z_{2}|\\
&\leq
|b|+\frac{1}{\zeta_i^2}(|z_{1}^{T}D_{1}^{-1}z_{1}|+|z_{2}^{T}D_{2}^{-1}z_{2}|),
\end{align*}
from (\ref{eq:508}) it follows
\begin{equation}
|\kappa _{\nu }|\leq 3 \sqrt{n}+(n+3)\bigg(\frac{|b|}{|\nu|}
+\frac{2}{|\nu|}\frac{\left\vert z_{1}^{T}D_{1}^{-1}z_{1}\right\vert
+\left\vert z_{2}^{T}D_{2}^{-1}z_{2}\right\vert) }{\zeta _{i}^{2}}\bigg) .
\label{eq:509}
\end{equation}%

Since $|b|\leq |\nu|$ and
\begin{align*}
\left\Vert A_{i}^{-1}\right\Vert _{2}& =|\nu|=\underset{\left\Vert x\right\Vert
_{2}=1}{\max }\left\Vert A_{i}^{-1}x\right\Vert _{2}\geq \left\Vert
A_{i}^{-1}e_{k}\right\Vert _{2} \\
& =\sqrt{\frac{1}{\left( d_{k}-d_{i}\right) ^{2}}+\frac{\zeta _{k}^{2}}{%
\zeta _{i}^{2}\left( d_{k}-d_{i}\right) ^{2}}}\geq \frac{\left\vert \zeta
_{k}\right\vert }{\left\vert \zeta _{i}\right\vert \left\vert
d_{k}-d_{i}\right\vert },
\end{align*}
by simply dividing each term 
$$
\frac{\zeta _{k}^{2}}{\zeta
_{i}^{2}|d_{k}-d_{i}|}
$$ 
in (\ref{eq:509}) 
with the corresponding quotient 
$$
\frac{\left\vert \zeta _{k}\right\vert }{%
\left\vert \zeta _{i}\right\vert \left\vert d_{k}-d_{i}\right\vert },
$$
we obtain 
\begin{equation}\label{konacno2}
\left\vert \kappa _{\nu }\right\vert \leq 3\sqrt{n}+(n+3)\bigg(1+\frac{2}{%
\zeta _{i}}\sum\limits_{\substack{ k=1  \\ k\neq i}}^{n-1}\left\vert \zeta
_{k}\right\vert \bigg) .
\end{equation}
The bound (\ref{knu2}) now follows from (\ref{konacno1}) and (\ref{konacno2}).
$\blacksquare$

\subsection*{Proof of Theorem \ref{T7}.}

We first prove the bound (\ref{kapa_mi_2}).
Since $\widetilde{\nu}=fl(\widehat{\nu})$ is computed by bisection, from
(\ref{lambda_max}) we have
$$
\widetilde{\nu}=\widehat{\nu}(1+\kappa_{bis}\varepsilon_M).
$$
This and (\ref{kapa_nu}) imply
$$
\widetilde{\nu}=\nu(1+\kappa_\nu\varepsilon_M) (1+\kappa_{bis}\varepsilon_M).
$$
Since $\widehat{\mu}=fl (1/\widehat{\nu})$, the bound 
(\ref{kapa_mi_2}) follows by ignoring higher order terms.
The bound (\ref{kapa_lambda_2}) now follows by inserting 
(\ref{kapa_mi_2}) into Theorems \ref{T1} and \ref{C1}.
$\blacksquare$

\subsection*{Proof of Theorem \ref{T6}.}

Let the assumptions of the theorem hold. Let $b$ be computed in double of 
the working precision, $\varepsilon_M^2$, and then stored in the standard
precision. 
The standard floating-point error analysis with neglecting higher order terms
gives 
\begin{align*}
\frac{P\left( 1+\kappa_P\varepsilon _{M}^{2}\right) - Q\left(
1+\kappa_Q \varepsilon _{M}^{2}\right)}{\zeta_i^2} \left( 1+\kappa_1\varepsilon
_{M}^{2}\right) & = \frac{P-Q}{\zeta_i^2} \left( 1+\kappa_b\varepsilon_M\right)
\\ &\equiv b\left( 1+\kappa_b\varepsilon_M\right),
\end{align*}%
where $|\kappa_P|,|\kappa_Q|\leq (n+1)$ and $|\kappa_1|\leq 3$.
Solving the above equality for $\kappa_b$, neglecting higher order terms,  and taking absolute values gives
\begin{equation*}
\left\vert \kappa_{b}\right\vert \leq \frac{\left\vert P\right\vert
+\left\vert Q\right\vert }{\left\vert P-Q\right\vert }\left( n+4\right)
\varepsilon _{M}
\equiv K_b (n+4)\varepsilon_M.
\end{equation*}%
Since, by assumption, $K_b\leq O(1/\varepsilon_M)$, this implies
$$
|\kappa_b|\leq O(n),
$$
as desired. $\blacksquare$


\end{document}